# A hybrid metaheuristic to optimize electric first-mile feeder services with charging synchronization constraints and customer rejections

Tai-Yu Ma[1], Yumeng Fang[1], Richard D. Connors[2], Francesco Viti[2], Haruko Nakao[2]

[1] Luxembourg Institute of Socio-Economic Research (LISER), 11 Porte des Sciences, 4366 Esch-sur-Alzette, Luxembourg

[2] Department of Engineering, University of Luxembourg, Esch-sur-Alzette, Luxembourg

**Abstract**

This paper addresses the on-demand meeting-point-based feeder electric bus routing and charging scheduling problem under charging synchronization constraints. The problem considered exhibits the structure of the location routing problem, which is more difficult to solve than many electric vehicle routing problems with capacitated charging stations. We propose to model the problem using a mixed-integer linear programming approach based on a layered graph structure. An efficient hybrid metaheuristic solution algorithm is proposed. A mixture of random and greedy partial charging scheduling strategies is used to find feasible charging schedules under the synchronization constraints. The algorithm is tested on instances with up to 100 customers and 49 bus stops/meeting points. The results show that the proposed algorithm provides near-optimal solutions within less one minute on average compared with the best solutions found by a mixed-integer linear programming solver set with a 4-hour computation time limit. A case study on a larger sized case with 1000 customers and 111 meeting points shows the proposed method is applicable to real-world situations.

## 1. Introduction

The climate crisis has brought the transport sector into a new era: the need to significantly reduce $CO_2$ emissions motivated massive investments for adopting non-combustion-fueled vehicles, and in particular electric vehicles (EVs). This emerging tendency brings new challenges for demand responsive transport (DRT) services, since current EVs need to recharge several times a day due to limited battery capacity (Jenn, 2019). When operating a fleet of EVs, efficient charging management becomes a critical component of the overall system cost. While EV routing problems have been extensively studied in the past decades, most studies assume unlimited charging station capacity due to the difficulty of solving problems that include capacitated charging (Froger et al., 2021). In addition, to improve system efficiency, the concept of meeting points has been adopted in several real-world microtransit and ridesharing services[1] (Haglund et al., 2019; Ma et al., 2021). In meeting-point-based DRT systems, customers are picked up and dropped off at nearby street corners or predefined feasible public transport stops within a reasonable walking distance from their origins or destinations (Czioska et al., 2019). Transport operators benefit from reducing their operation costs with little increase of customer's inconvenience. Several recent studies have applied the meeting-point-based concept to improve the operational efficiency of DRT in both static and dynamic settings (Melis and Sörensen, 2022; Montenegro et al., 2022). However, integrating electric vehicles into this kind of on-demand DRT system has not yet been studied.

Deploying electric vehicles in a DRT system needs to address the joint optimization problems of vehicle routing and charging scheduling. This problem is closely related to the electric door-to-door DRT system (Pimenta et al., 2017; Bongiovanni et al., 2019). The problem consists of deciding vehicle routes and charging schedules to serve a set of customers while satisfying constraints on vehicle capacity, time windows, and vehicle energy consumption (Kucukoglu et al., 2021). Most of the literature assumes that vehicles can be charged at any time with unlimited charging station capacity (Schneider et al., 2014). This assumption is often violated in practice, as the number of fast chargers is very limited due to their high installation costs. The Electric Vehicle Routing Problem with Capacitated Charging Stations (EVRP-CS) is even more challenging as it needs to synchronize the charging operations of vehicles to save waiting time at charging stations. Recent research efforts have mainly focused on developing exact methods based on mixed linear integer programming (MILP) by assuming that vehicles are fully charged after leaving charging stations (Brugliei et al., 2019; Froger et al., 2021).

---

[1] https://blog.blablacar.com/blog/blablalife/travel-tips/ridesharing-meeting-points

In this paper, we focus on a flexible meeting-point-based electric DRT (feeder) system which provides a passenger transport service to connect to transit stations. This type of service is mainly applied in rural areas where public transport service is poor (Alonso-González et al., 2018; Ma et al., 2021). Moreover, customers may be rejected to consider the tradeoffs of operation costs and level of service. The problem needs to decide jointly where to pick up customers and how to route vehicles. The problem is complex due to interactions between customer-to-meeting-point assignment and the subsequent vehicle routing and charging synchronization under charging station capacity constraints.

The contributions of this paper are summarized as follows.
- We propose a MILP model to address the meeting-point-based electric feeder service problem with charging synchronization constraints (MP-EFCS) and allowing customer rejections. This problem extends the existing electric dial-a-ride problem (e-DARP) (Bongiovanni et al., 2019) by allowing multiple partial recharges of vehicles and including charging station capacity constraints. The latter is modeled at charger level, as a job scheduling constraint to ensure no vehicle charging conflicts occur.
- To solve the MILP model efficiently, we propose a new layered graph structure to trim unnecessary nodes and arcs to reduce the problem size. This layered graph structure is similar to the time-expanded graph, but we use the information of their layer index to help reduce the problem size for the customer to meeting point assignment.
- A hybrid metaheuristic algorithm is proposed to solve the MP-EFCS problem efficiently. This algorithm first assigns customers to meeting points, as a variant of the capacitated facility location problem, and then applies a deterministic-annealing-based (DA) algorithm to solve our variant of the EVRPTW-CS problem. The resulting solution is further optimized by a matheuristic for customer reassignment when the solution contains unserved customers.
- Two sets of benchmark instances with up to 100 customers are generated to evaluate the performance of the algorithms and compare them with solutions obtained by a commercial solver. Different initial battery levels and demand distributions are considered providing a more general evaluation of the performance of the proposed solution method.
- To ensure good performance, we conduct a sensitivity analysis to set up the algorithmic parameters. A case study is conducted to analyze the impact of system parameters with respect to the meeting point separation distance and fleet size on the system performance.
- Finally, the experiment is extended for larger instances corresponding to a real-world case in the Arlon-Luxembourg cross-border area with 1000 requests. We show that our algorithm is applicable for solving real-world larger instances in reasonable time. Both system performance and computational results are analyzed.

The remainder of the paper is organized as follows. Section 2 reviews the related literature and highlights the research gaps. Section 3 presents the problem description and its MILP formulation. Section 4 presents the hybrid metaheuristic algorithm, which consists of three subproblems: 1) customer-to-meeting-point assignment; 2) e-DARP with charging synchronization (e-DARP-CS); 3) customer reinsertion. A layered graph model is proposed to prune the problem size; an illustrated example illustrates the method. Section 5 presents the computational study, including test instance generation, algorithm parameter settings, computational results, and case study. Finally, conclusions are drawn and future extensions discussed.

## 2. Related studies

In this section, we review the related literature focusing on the following aspects: Demand-responsive feeder service, meeting-point-based models, and electric vehicle routing with charging synchronization constraints. The reader is referred to Vansteenwegen et al. (2022) for a comprehensive review of DRT systems. Research gaps are summarized at the end of this section.

### 2.1. Demand-responsive feeder service

Demand-responsive feeder systems provide on-demand dial-a-ride service to increase public transport ridership (e.g. by transporting customers to transit stations). It is a cost-effective solution in rural areas where public transportation services are not well developed. Lee and Savelsbergh (2017) reviewed previous works and classified the full spectrum of DRT services according to different flexible route design concepts (zone-based vs. point-based deviation, route deviation, flexible route segments, or demand-responsive connection to a transit hub). The authors formulate the demand-responsive connector problem as a MILP to minimize the total routing cost. When booking a trip, each customer indicates their earliest pickup time and specifies their latest arrival time at the transit station (corresponding to scheduled transit departure time). The transit connector service ensures that there are no late drop-offs so that customers do not miss their trains. In addition, customers can be dropped off at alternative stations instead of predetermined transit stations (regional systems). Their computational results show that this flexibility provides cost savings of up to about 29% compared to the conventional regional systems. Chen and Nie (2017) analyze an idealized demand-responsive connector system that uses demand-adaptive services to connect to fixed-route transit lines for customer transfers. They propose an analytical model to evaluate the impact of various system design parameters on system performance, including e.g. road spacing distance, value of travel time, walking distance, vehicle speed, and operating cost per vehicle mile/hour. Montenegro et al. (2022) further consider two types of bus stops: mandatory stops (must be visited by buses) and optional stops (visited upon there are customers nearby). Customers are assigned to the optional stops within a predefined maximum walking distance. The objective is to minimize a weighted sum of bus routing times, customer walking times, and the penalty associated with early/late arrivals with respect to the customer's desired arrival times. Bian and Liu (2019) analyze the mechanism of first-mile ridesharing service from both the operator (customer-route matching) and customer (incentives, customized pricing) perspectives. Several studies evaluate the performance of demand-responsive feeder service using simulation approaches or empirical trip data from operators (see, for example, Alonso-González et al., 2018; Ma et al., 2019; Haglund et al., 2019; Wang and Ross, 2019, among many others). However, existing studies assume that all requests need to be served with sufficient fleet size. This assumption might not be practical for operators when operation costs for serving certain customers are relatively high. Few studies allow customers to be rejected by considering the trade-off in the objective function.

### 2.2. Meeting-point-based models for DRT system design

To improve the efficiency of door-to-door-based DRT systems, meeting-point-based approaches have recently received increasing interest and have been applied in real-world ridesharing/microtransit services (Aïvodji et al., 2016; Stiglic et al., 2015; Zheng et al., 2019). These approaches consist of picking up and dropping off customers at nearby, predefined "meeting points," which are places/areas that make it easy for customers to get on and off vehicles. Czioska et al. (2019) reviewed previous work on the methods used to define the meeting points, namely using a set of predefined feasible locations or dynamically generating them based on the locations of requests. Once the customer meeting points are assigned, the operator optimizes the vehicle routing to minimize the total routing cost. The authors propose a three-step procedure by first clustering requests based on their origin-destination pairs, then generating meeting points within a set of feasible locations (parking lots, side street intersections, gas stations, etc.). A meeting point is selected by checking whether it is feasible for all customers on that trip. Finally, vehicle routes are optimized by solving a vehicle routing problem based on the selected meeting points. However, this sequential approach leads to suboptimal solutions because these decisions are interdependent. This problem is a variant of school bus routing problems (Park and Kim, 2010; Schittekat et al., 2013), but is more challenging due to the need to satisfy time windows and travel time constraints of customers. Given its potential utility for DRT systems, Montenegro et al. (2021, 2022) propose a demand-responsive feeder service with both mandatory and optional bus stops (meeting points), where the assignment of customers to optional bus stops and vehicle routing are jointly optimized. The authors develop a column generation method to find optimal solutions for small instances and a large neighborhood search approach to solve large instances. Melis et al. (Melis et al., 2021; Melis and Sörensen, 2022) address a similar problem called on-demand bus routing problems to jointly optimize customer-to-meeting-point assignment and bus routing to minimize the total travel time of users. The authors point out that such a system has good potential to reduce system costs, especially in the context of autonomous vehicles. However, these studies assume that all customers need to be served and are based on a fleet of homogeneous internal combustion engine vehicles. In the context of electric vehicles, the problem becomes more complex where vehicle energy constraint needs to be satisfied and charging schedules need to be jointly optimized.

2.3. Electric vehicle routing problems with capacitated charging stations

Electric vehicle routing problems extend classical vehicle routing problems by considering vehicle battery range limits and constraints on charging operations to minimize overall routing and charging costs (Erdoğan and Miller-Hooks, 2012). Since vehicle range is limited, vehicles need to be recharged en route while minimizing the impact on customer service. Previous studies range from assuming linear charging functions and full recharging (Schneider et al., 2014) to more realistic nonlinear charging functions and allowing partial recharging (Desaulniers et al., 2016). For example, Felipe et al. (2014) consider a partial-recharge policy where the charging amount depends on the remaining trips of vehicles to minimize the charging time and cost. This policy results in lower operating costs and greater vehicle availability to serve more customers. Similar conclusions have been drawn regarding the benefits of adopting a partial recharge policy (Keskin and Çatay, 2016). Regarding the modeling of the charging function, most studies assume a linear function with constant charging rate. However, a more realistic charging behavior is nonlinear, the charging speed is significantly slowed down when the battery is above 80% of its capacity. To account for this aspect, several studies approximate the nonlinear charging speeds with piecewise linear functions to obtain more accurate charging times and costs when modeling related problems (Keskin et al., 2019; Froger et al., 2021; Lam et al., 2022). From this perspective, methods for realistic energy consumption function estimation incorporating vehicle speed, road profile, vehicle load, and acceleration/deceleration, are also studied (Goeke and Schneider, 2015; Macrina et al., 2019).

While electric vehicle routing problems have been studied extensively, only a few studies consider capacity constrained charging stations. Bruglieri et al. (2019) propose two different modeling approaches (arc-based and path-based) to formulate the charging station capacity constraints in electric vehicle routing problems. To allow multiple visits of vehicles to charging stations, multiple dummy copies are created for each charging plug. Charging capacity is ensured by postponing the current visit of a charger by at least the full charging time of the previous visit of another vehicle in ascending order of visits. The authors propose an alternative path-based MILP formulation and use the cutting plane method to solve the test case with less than 20 customers to the optimal solution in less than 1 hour of computation time (Bruglieri et al., 2021), Froger et al. (2021) propose a path-based MILP formulation by considering piecewise linear charging functions and partial recharging, and propose a matheuristic method to solve the EVRP-CS exactly. Their solution method first generates a pool of initial routes without considering the capacity constraints of the charging stations, and then in the second step, they try to recombine these routes to find a solution that satisfies the capacity constraints. Their problem assumes that the vehicles are homogeneous (in terms of battery size) and fully charged before starting the service. They were able to solve most of the test instances with 10 customers exactly. Lam et al. (2022) propose a branch-and-cut-and-price algorithm to solve the EVRP problem with time windows and capacitated charging station constraints (EVRPTW-CS) by considering both partial vehicle recharge and piecewise linear charging functions. The subproblem of synchronization of charging schedules is solved by applying the constraint programming technique. Its exact method can solve the problem with some larger test instances of up to 100 customers. As for the heuristic approaches, Keskin et al. (2019) develop an adaptive large neighborhood search algorithm to solve EVRP considering the waiting time at charging stations for the problem instances of up to 100 customers. To the best of our knowledge, the existing literature mainly focuses on exact methods that can be applied to small instances. There are still no efficient algorithms for large-scale electrical dial-a-ride problems with capacitated charging stations.

In summary, an on-demand, meeting-point-based bus system has good potential to achieve significant system cost savings. However, existing studies have not yet considered an electric fleet, for which vehicle charging synchronization remains a challenging issue. This study aims to address these issues to develop a solution approach for this type of system in the context of on-demand feeder service using a fleet of heterogeneous electric vehicles.

## 3. Problem description and formulation

We consider a DRT feeder service operated by an operator in a rural area using a heterogeneous (in terms of capacity, battery size, and energy consumption rate) fleet of electric buses (also called vehicles hereafter) to complement the public transport system. To enhance system efficiency and reduce operational costs, the DRT system adopts the concept of **meeting points** i.e. customers are offered a limited number of pick-up/drop-off meeting points, rather than a door-to-door service (Czioska et al., 2019; Ma et al., 2021) and the service is **punctuated** (e.g. the vehicle arrives at a transit station every 10-20 minutes to drop off the transit passengers). The system is operated as follows. For a given planning period, customers submit their ride requests in advance indicating their origin, the transit station to be dropped off, and their desired arrival time (corresponding to the pre-defined departure time of the transit service). Each request (customer) contains at least one passenger. The operator collects these ride requests and communicates whether they are accepted, the pickup time, and suggested meeting points. The operator's objective is to optimize vehicle routes so as to arrive at transit stations within a fixed buffer time (e.g. $\leq 10$ minutes before the timetabled transit departure). We assume that customers are willing to walk from their origins to the suggested meeting points, up to some maximum acceptable walking distance. The state of charge of the vehicles cannot fall below the reserve battery level throughout the route. Vehicles can be recharged only at operator-owned charging stations; each station has a limited number of chargers. Charging operations cannot overlap at any charger i.e. waiting of a vehicle is not allowed at a charger/charging station. Given a set of customer requests, the objective is to optimize vehicle routes to meet these requests while considering the trade-off between system costs and customer inconvenience.

The MP-EFCS problem is formulated as a mixed integer linear programming (MILP) problem. Note that 3-index formulation is necessary since we consider a fleet of heterogeneous vehicles with different initial states of charge and vehicle energy consumption rate. For simplifying the analysis, we formulate our model in a single-depot setting. This can be extended to a multi-depot model without difficulty by indicating the starting and returning depots of each vehicle (Bongiovanni et al., 2019; Braekers et al., 2014). The objective function minimizes the weighted sum of total vehicle travel time and total vehicle charging time (the first term), customer's total walking time (the second term), total 'excess' vehicle waiting time at transit stations (exceeding the acceptable fixed buffer time), and the total penalty for unserved requests. The weighting factors $\lambda_1, \lambda_2, \lambda_3, \lambda_4$ are user-specified parameters to account for trade-offs between users' and operators' objectives. The notation table can be found in Appendix A.

$$\text{Min } Z = \lambda_1 \sum_{k \in K} \left( \sum_{(i,j) \in \mathcal{A}_B} t_{ij} x_{ij}^k + \sum_{s \in S'} \tau_s^k \right) + \lambda_2 \sum_{(r,i) \in \mathcal{A}_c} \sum_{k \in K} y_{ri}^k t_{ri}$$
$$+ \lambda_3 \sum_{k \in K} \sum_{i \in D'} W_i^k + \omega \sum_{(r,i) \in \mathcal{A}_c} \left(1 - \sum_{k \in K} y_{ri}^k\right) \quad (1)$$

The constraints can be grouped into four categories in terms of customer-to-meeting-point assignment constraints, vehicle routing constraints, vehicle energy constraints, and charging scheduling (synchronization) constraints. Equation (2) ensures that each customer is served at most once. Equation (3) imposes the maximum walking distance of customers to access a meeting point.

$$\sum_{k \in K} \sum_{i \in G'} y_{ri}^k \leq 1, \quad \forall r \in R \quad (2)$$

$$\sum_{k \in K} \sum_{i \in G'} w_{ri} y_{ri}^k \leq w_{max}, \quad \forall r \in R \quad (3)$$

In terms of vehicle routing constraints, equations (4) and (5) state that each vehicle leaves the depot and returns to the same depot. Equations (6) ensures that each meeting point (dummy) node can be visited at most once by the same vehicle. Equation (7) ensures vehicle flow conservation. Equation (8) ensures consistency between $y_{ri}^k$ and $x_{ij}^k$. Equation (9) ensures that the pickup and drop-off of a customer is served by the same vehicle. Equations (10)-(11) update the bus occupancy at meeting points (pick-up locations) and transit stations (drop-off locations). Equation (12) states the capacity (passenger load) constraint of the vehicle. Equation (13) states that the beginning time of service at node $j$ can start when the bus arrives at $j$. Equation (14) states that when leaving a charger $s$, the starting time of service at successive node $j$ is constrained by its starting time of service at the preceding node s plus the service time, charging time, and the travel time traversing arc $(i, j)$. Equation (15) computes the arrival time of vehicles at transit stations. Note that the hard time windows (fixed buffer time associated with the transit service timetable) are associated with transit station nodes only, not for meeting point nodes. To determine the excess waiting time when arriving at a transit station (before the buffer time), we introduce the arrival time variable, relevant only for the transit stations. Equation (16) then measures the excess bus waiting time. Equation (17) determines the value of the auxiliary variable indicating whether there are buses dropping off customers at transit station node $i$. Equation (18) ensures the ride time of customers cannot exceed the maximum ride time, characterized as the shortest travel time (direct ride) multiplied by a pre-defined detour factor. Equation (19) states that the beginning time of service at node $i$ is constrained by the time window associated with that node.

$$\sum_{j \in G' \cup S' \cup \{N+1\}} x_{0j}^k = 1, \forall k \in K \quad (4)$$

$$\sum_{i \in \{0\} \cup S' \cup D'} x_{i,N+1}^k = 1, \forall k \in K \quad (5)$$

$$\sum_{i \in V_0} x_{ij}^k \leq 1, \quad \forall k \in K, j \in G' \quad (6)$$

$$\sum_{i \in V_0} x_{ij}^k - \sum_{i \in V_{N+1}} x_{ji}^k = 0, \quad \forall k \in K, j \in V \quad (7)$$

$$\sum_{r \in R} y_{ri}^k \leq M_1 \sum_{j \in V_{N+1}} x_{ij}^k, \quad \forall k \in K, i \in G' \quad (8)$$

$$y_{ri}^k = 1 \Rightarrow \sum_{j \in V_0} x_{ji}^k = \sum_{j \in G' \cup D'} x_{jd_r}^k, \quad \forall k \in K, i \in G', r \in R \quad (9)$$

$$x_{ij}^k = 1 \Rightarrow q_j^k = q_i^k + \sum_{r \in R} y_{rj}^k, \quad \forall k \in K, i \in V_0, j \in G' \quad (10)$$

$$x_{ij}^k = 1 \Rightarrow q_j^k = q_i^k - \sum_{r \in R} \sum_{g \in G'} y_{rg}^k, \quad \forall k \in K, i \in G', j \in D' \tag{11}$$

$$0 \leq q_i^k \leq Q^k, \quad \forall k \in K, i \in V_{0,N+1} \tag{12}$$

$$B_j^k \geq B_i^k + u_i + t_{ij} - M_2(1 - x_{ij}^k), \quad \forall k \in K, i \in V_0, j \in V_{N+1} \tag{13}$$

$$B_j^k \geq B_s^k + \tau_s^k + t_{sj} - M_2(1 - x_{sj}^k), \quad \forall k \in K, s \in S', j \in \{G' \cup N+1\} \tag{14}$$

$$x_{ij}^k = 1 \Rightarrow A_j^k = B_i^k + t_{ij} + u_i, \quad \forall k \in K, i \in G' \cup D', j \in D' \tag{15}$$

$$W_i^k \geq B_i^k - A_i^k - M_2(1 - p_i^k), \quad \forall k \in K, i \in D' \tag{16}$$

$$p_i^k = \sum_{j \in V} x_{ji}^k, \quad \forall k \in K, i \in D' \tag{17}$$

$$A_{d_r}^k - B_i^k - u_i \leq L_i + M_2(1 - y_{ri}^k), \quad \forall k \in K, (r,i) \in A_C \tag{18}$$

$$e_i \leq B_i^k \leq l_i, \quad \forall k \in K, i \in V \tag{19}$$

For vehicle energy constraints, equations (20) and (21) state the initial battery level of the vehicles and their constraints. Equations (22)-(23) ensure energy conservation with and without recharged energy when traversing an arc $(i, j)$.

$$E_0^k = E_{init}^k, \quad \forall k \in K \tag{20}$$

$$E_{min}^k \leq E_i^k \leq E_{max}^k, \quad \forall k \in K, i \in V \tag{21}$$

$$x_{ij}^k = 1 \Rightarrow E_j^k = E_i^k - \beta^k c_{ij}, \quad \forall k \in K, i \in V_0 \backslash S', j \in V_{N+1} \tag{22}$$

$$x_{ij}^k = 1 \Rightarrow E_j^k = E_s^k + \alpha_s \tau_s^k - \beta^k c_{sj}, \quad \forall k \in K, s \in S', j \in \{G' \cup N+1\} \tag{23}$$

With respect to the constraints on charging scheduling (synchronization), equations (24)-(28) ensure that each charger can be occupied by no more than one vehicle at a time, i.e., if multiple charging events are scheduled at the same charger, they cannot overlap. Creation of dummy charger nodes allows multiple visits by vehicles to the same chargers. Equation (24) introduces an auxiliary variable to indicate whether a dummy charger node is visited or not. Equation (25) states that the dummy charger nodes are visited in reverse order as they appear in the list of their associated physical charger nodes to eliminate the symmetry problem (see Figure 1) (see e.g. Froger et al., 2017; Lee and Savelsbergh, 2017). Equation (26) ensures that charging of a vehicle can start only after the previous charging has finished. This means that the start time of a charging visit cannot be earlier than the start time of the previous charging visit of a vehicle plus its charging duration. In equation (27), if a vehicle is not connected to a dummy charger node, its start time and charge duration are set to zero to determine which vehicle is connected to the dummy charger node in equation (26). Since each dummy charger node can only be visited once, the unique connected vehicle and its associated charge start time and charge duration can be determined by (26). Equation (28) states that each dummy charger node can be visited at most once. Note that if the symmetry issue is not addressed when modeling multiple visits to charging stations, the computation time will increase significantly even for moderate problem sizes.

$$v_s = \sum_{k \in K} \sum_{j \in G' \cup N+1} x_{sj}^k, s \in S' \tag{24}$$

$$v_h \leq v_l, \forall h, l \in S'_o, o \in S, h < l \tag{25}$$

$$\sum_{k \in K} B_h^k \geq \sum_{k \in K} B_l^k + \sum_{k \in K} \tau_l^k - M_2(2 - v_h - v_l), \forall h, l \in S'_o, o \in S, h < l \tag{26}$$

$$\tau_s^k + B_s^k \leq M_2 \sum_{j \in G' \cup N+1} x_{sj}^k, \forall s \in S', k \in K \tag{27}$$

$$v_s \leq 1, s \in S' \tag{28}$$

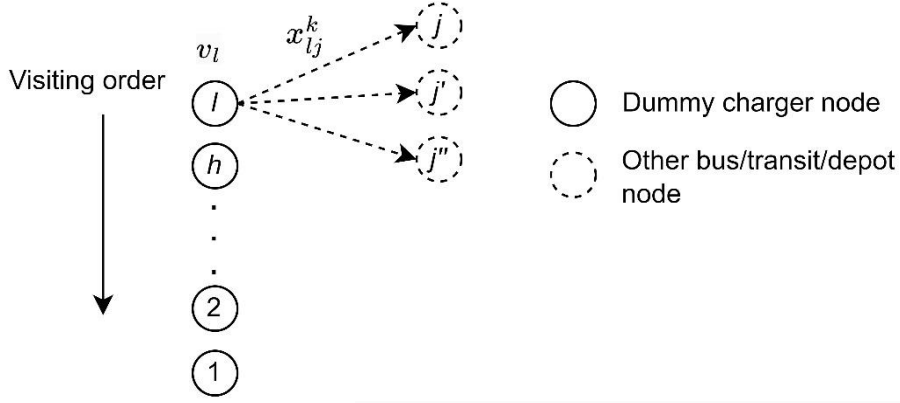

**Figure 1. Illustration of the visit order of dummy charger nodes associated with a physical charger. Each dummy node can be connected maximum once by vehicles, and follows the inverse order (from node *l*, then node *h*, (*i.e., l*-1), ..., 2,1.**

Finally, equations (29)-(33) define the domain of the decision/auxiliary variables.

$$x_{ij}^k \in \{0,1\}, \quad \forall k \in K, i,j \in V_{0,N+1} \tag{29}$$
$$y_{ri}^k \in \{0,1\}, \quad \forall k \in K, r \in R, i \in G' \tag{30}$$
$$\tau_s^k \geq 0, v_s \in \{0,1\}, \forall k \in K, s \in S' \tag{31}$$
$$A_i^k \geq 0, B_i^k \geq 0, \forall k \in K, i \in V_{0,N+1} \tag{32}$$
$$p_i^k \in \{0,1\}, W_i^k \geq 0, \forall k \in K, i \in D' \tag{33}$$

Note that Equation (9) can be re-written equivalently as the following two equations.

$$\sum_{j \in V_0} x_{ji}^k \geq \sum_{j \in G' \cup D'} x_{jd^r}^k - M_1(1 - y_{ri}^k), \quad \forall k \in K, i \in G', r \in R \tag{34}$$

$$\sum_{j \in V_0} x_{ji}^k \leq \sum_{j \in G' \cup D'} x_{jd^r}^k + M_1(1 - y_{ri}^k), \quad \forall k \in K, i \in G', r \in R \tag{35}$$

Eqs. (10)-(11), (15), and (22)-(23) can be extended equivalently in the similar way as above. The big positive numbers $M$ are set as $M_1 = |R|, M_2 = T$.

Obviously, the number of dummy nodes for each meeting point (same location) need to be at least equal to the number of potential visits, corresponding to the requested bus service arrivals, at transit stations. To solve this problem, the dummy nodes need be sorted according to a desired structure to reduce the search space of the problem.

### 4. Hybrid metaheuristic algorithm

The MP-EFCS problem described in the previous section has structure similar to the location-routing problem (Belenguer et al., 2011) or on-demand bus routing problems (Melis et al., 2021; Montenegro et al., 2022) in which customers need to be assigned to nearby meeting points first, and then bus routes are optimized under vehicle capacity, battery level, customer ride time, drop-off time windows and charging synchronization constraints. We propose an efficient **two-stage solution scheme** by finding a good customer-to-meeting-point assignment in the first stage. In the second stage, a metaheuristic with a post-optimization procedure is proposed to solve the routing problem with charging synchronization constraints. The proposed metaheuristic is hybrid because it incorporates two additional subproblems (customer-to-meeting-point assignment problem and post-optimization problem, each formulated as a MILP) on top of the main EVRP-CS routing problem into the overall solution framework. The goal is to avoid solving the EVRP-CS problem multiple times, each corresponding to a possible customer-to-meeting-point assignment outcome. Notice that different from existing e-DARP problem formulations (see e.g.

Bongiovanni et al., 2019) and on-demand bus routing problems (Montenegro et al., 2022), the proposed model is more flexible since vehicles can be partially recharged, vehicles' waiting times at transit stations need to be minimized and customers can be rejected if all constraints cannot be met. We propose an DA-based metaheuristic that has been successfully used in solving the (gasoline vehicle) general heterogeneous dial-a-ride problem (Braekers et al., 2014). The proposed DA algorithm is similar in structure to Braekers et al. (2014), but tailored to the problem of MP-EFCS with additional procedures to reduce the computation time. The DA applies a number of local search operators in the neighborhood of the current solution to improve it. A temperature parameter is used as a threshold to allow accepting worse solutions so as to escape from local optima. This temperature parameter is gradually reduced during the optimization process until no further improvement, or some stopping criteria, are met. Compared to other metaheuristics, such as variable neighborhood search (Schneider et al., 2014), the advantage of the DA is that there are only a small number of algorithmic parameters to be tuned. Furthermore, the DA performs nearly equally well as simulated annealing algorithm but needs only a fractional computational time compared to the latter (Dueck and Scheuer, 1990). The new challenge is how to optimize vehicle charging schedules with synchronization constraints, and how to improve the quality of the solution by efficiently reassigning customers to different meeting points when certain customers are unserved. This will be explained in the following sections and in the computational studies.

4.1. Preprocessing

The MP-EFCS problem in Section 3 is defined on a directed graph where the set of nodes includes the set of customers R and the set of meeting point nodes ($G'$), transit station nodes ($D'$), charger nodes ($S'$) and the two copies of the depot {0, N+1}. There are two sets of arcs: 1) walking arcs $\mathcal{A}_c$ for customers walking from their origins to the meeting point nodes, and bus (vehicle) arcs $\mathcal{A}_B$ for bus routes starting from the depot node 0 and terminating at $N + 1$. Each arc is associated with a distance and a travel time, calculated as the Euclidean distance divided by the average walking/vehicle speed. Different from the classical dial-a-ride, where customer's pick-up and drop-off locations are given, customers are first assigned to nearby meeting points, following which, bus routes and charging scheduling can be optimized. The operator then communicates the pick-up times at meeting points (stops of the buses) and the assigned meeting points to the customers, if their request for a ride has been accepted. It is assumed that customers have real-time information about the arrival time of the buses, so the waiting time of customers can be ignored. Time windows are associated with transit stations with the aim of obtaining solutions with low waiting time for the transfer of customers. Solving the MP-EFCS exactly necessitates enumerating all possible customer-to-meeting-point assignments, and then solving each corresponding e-DARP-CS problem to find the global minimum. This is possible only for very small problem sizes. To solve it efficiently, we propose a **layered (directed) graph** model according to the sorted arrival timetable at transit stations, and prune infeasible arcs and unnecessary nodes to reduce the problem size. The following definitions are used for the layered graph.

- **Layer**. A layer is a subset of meeting point and transit station dummy nodes generated for each pair of transit station and ordered index of arrival time of the service timetable (see an example in Figure 2). The index set of the layers is defined as $\ell \in \{1,2,\ldots,\mathcal{L}\}$, where the latest arrival times at the higher layers are not less than those at the lower layers. A particular layer is characterized by the pair: scheduled arrival time and associated transit station.
- **Compatible layers**: Two layers are compatible if the meeting point nodes (or transit station nodes) on the two layers can be visited by the same vehicle, i.e. given a meeting point, there is a vertical arc connecting a dummy node on one level to the one on the other level and vice versa. This compatibility can be verified by checking whether the arcs connecting the two dummy transit station nodes of the two layers can be part of the feasible solution. Let $[e_i, l_i]$ and $[e_j, l_j]$ denote the beginning and end of the time window at the dummy transit station node $i$ and $j$, respectively. The width of the time window corresponds to the predefined buffer time. Two layers $(i,j), i < j$ are compatible if $e_i + t_{ij} \leq l_j$. If two layers are compatible, there are vertical arcs connecting the dummy nodes of the same meeting point or same transit station.

- **Layered graph**: We denote a layered graph as $\mathcal{G} = (V_{0,N+1}, \mathcal{A})$, where $V_{0,N+1}$ is a set of nodes structured with a ground layer with two copies of the depot and the dummy charger nodes, and a set of the layers with meeting point and transit station dummy nodes (see the definition in Appendix A), sorted according to the increasing order of latest arrival time at transit station dummies $l_j, \forall j \in D'$ (dummy transit nodes). $\mathcal{A} = \mathcal{A}_B \cup \mathcal{A}_C$ is a set of arcs where $\mathcal{A}_B$ is a set of arcs for bus routing after trimming infeasible arcs, and $\mathcal{A}_C$ is a set of (walking) arcs connecting customer's origins to reachable meeting points within a predefined maximum walking distance. Note that $\mathcal{A}_B$ contains the following arcs:
  - Arcs from node 0 to $G'$, $S'$, and $N+1$
  - Arcs from $D'$ to $S'$ and $N+1$
  - Arcs connect from meeting point nodes to different meeting point nodes and transit station nodes on the same layer
  - Arcs from transit station nodes to meeting point nodes on the higher layers
  - Arcs connect the meeting point (transit station) nodes of the same physical location on two different layers if their layers are compatible
  - Arcs from $S'$ to $G' \cup N+1$
  
  $\mathcal{A}_C$ contains the arcs within the maximum walking distance of customers on the layer of the customer's desired train departures (see Figure 2).

Although the problem size is significantly reduced, the considered problem is NP-hard and can be solved exactly only with tens of customers (requests). Figure 2 shows an example of a layered graph for the e-DARP-CS problem. Depending on initial battery levels, some buses might need to recharge after dropping off customers at transit stations. The duplicate meeting point nodes allow the same vehicle to visit the same physical meeting point at different times. Layers without customer requests are trimmed. Let $\ell \in \tilde{\mathcal{L}}$ be the subset of layers with positive customer requests and $V_{\tilde{\mathcal{L}}} = G'_{\tilde{\mathcal{L}}} \cup D'_{\tilde{\mathcal{L}}}$ be the subset of "active" dummy nodes on these layers. Given any $\ell \in \tilde{\mathcal{L}}$, the active dummy meeting point nodes concern only those within the maximum customer walking distance ($w_{max}$) for each $\ell \in \tilde{\mathcal{L}}$. Thus, the problem size of MP-EFCS is a function of the number of active nodes which depend on the parameters of $w_{max}$, number of meeting points based on the methods to generate them, and bus service frequency to connect transit stations.

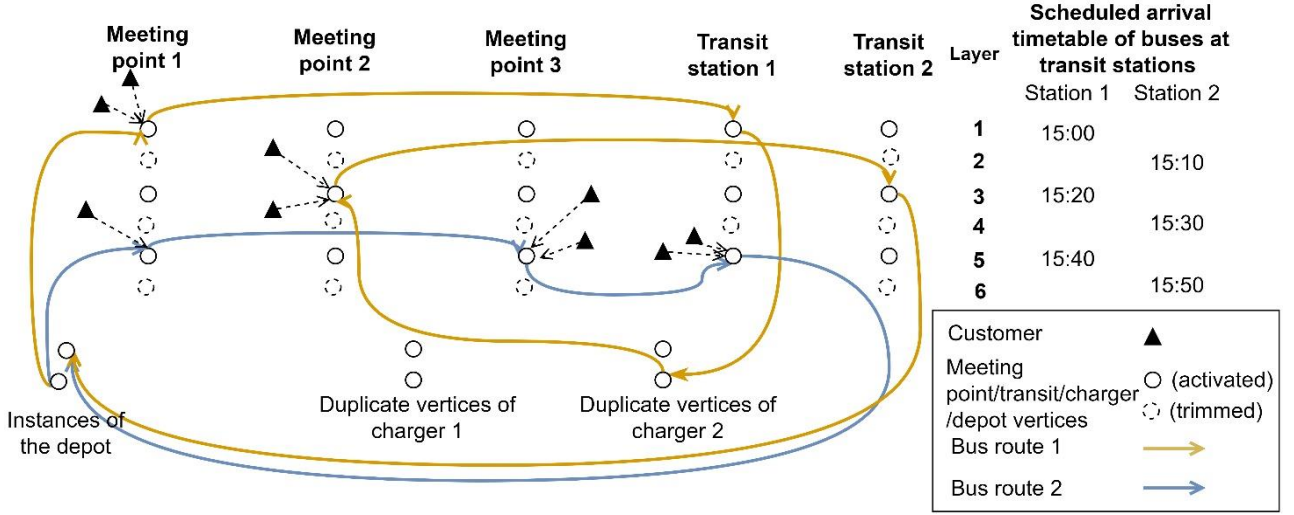

**Figure 2. An illustrative example of the layered directed graph (arcs are omitted) for modeling the meeting-points-based electric feeder service with the charging synchronization constraints.**

4.2. Hybrid metaheuristic algorithm

The overall structure of the algorithm is shown in Algorithm 1. It consists of three parts: customer-to-meeting-point assignment, electric vehicle routing, and post-optimization to insert unserved customers.

Starting from the input of data and parameters, a layered graph of the problem instance is constructed by trimming unnecessary/infeasible nodes/arcs. First, a customer-to-meeting-point assignment is made (line 2 in Algorithm 1, described later) then an e-DARP-CS instance is constructed by trimming off unused meeting point nodes and arcs connected to them, based on the layered graph model. Let $s$ denote the current solution, and $s_{best}$ denote the current best solution. $c(s)$ is the cost (objective function value of Equation (1)) of solution $s$. An initial feasible solution $s_{init}$ is generated as the best feasible solution found for $n$ random solutions (i.e. n=100) using a greedy insertion approach (line 4 in Algorithm 1). If this fails, $n$ is increased to 1000 to find a feasible solution. The algorithm applies a randomly selected local search operator $ls$ on the current solution $s$ and obtains a temporary solution $s'$ (lines 11-12 in Algorithm 1). If the cost of $s'$ is smaller than that of the current solution $s$ plus a threshold value $T$ (i.e. a worse solution but within the threshold), and there are no charging operation conflicts, a vehicle exchange operator (line 14 in Algorithm 1) is applied on $s'$ to further reduce the charging time of the vehicles of $s'$. This is because we assume that the initial battery levels of the vehicles are heterogeneous. Given two routes, this route switcher could reduce the total charging time if the vehicle with the higher battery level travels the longer route, thus reducing the amount of energy to be charged. To handle this operation efficiently, this operator sorts the vehicles by their total charging time, including the additional time to access the charging stations. A first vehicle is exchanged with a second vehicle (without changing the sequence of pickups and drop-offs of the route). If the resulting vehicle exchange and rescheduled charging operations (if any) improves cost without charging conflicts for all vehicles, then the current solution is updated. When the cost of $s$ is smaller than that of the best solution $s_{best}$, and the number of used vehicles $n_k(s)$ does not exceed the fleet size $|K|$, update the best solution (lines 16-18 in Algorithm 1), and reset the non-improvement count $i_{imp} = 0$. When $s_{best}$ is not improved ($i_{imp} > 0$), reduce the threshold value by $T_{max}/T_{red}$, where $T_{red}$ is the threshold reduction parameter. If $T < 0$, reset $T$ randomly between 0 and $T_{max}$. When $i_{imp}$ exceeds the maximum number of iterations, reset $i_{imp} = 0$ and the current solution as $s_{best}$.

We track the number of times (a check is applied once for every 100 iterations) that the current best solution has stagnated. If this exceeds a pre-defined limit, the algorithm returns the current best solution (lines 8-10 Algorithm 1). This is controlled by a user-defined parameter $n_{stagnant}$. Otherwise, randomly selected local search operators are applied until a maximum number of iterations is achieved. This early stop criterion helps to reduce the computation time. Note that $T_{max}, T_{red}, n_{imp}$ and $n_{stagnant}$ are the algorithmic parameters to be tuned. Note that the proposed algorithm can be adapted with little effort for multi-depot cases by specifying the depot locations of each vehicle.

The e-DARP-CS instance is optimized based on the solution obtained from the first stage customer-to-meeting-point assignment problem. It might be possible to accommodate unserved customers by changing their assigned meeting points, then reinserting them into the current bus routes. In doing so, bus routes and charging schedules need to be updated accordingly. In the case that there are unserved customers, we propose an efficient post-optimization procedure to re-optimize the best solution obtained from the DA algorithm (line 30 Algorithm 1). Our numerical results show that this post-optimization procedure can improve the final solution and reduce the number of unserved customers with reasonable additional computational effort. In the following, we present the mathematical formulations and algorithmic description of the three main components of the hybrid metaheuristic.

**Algorithm 1. Hybrid metaheuristic algorithm for solving meeting-point-based electric feeder service problem.**

```
1. 𝒢 ← createLayeredGraph($V_{0,N+1}$)
2. $(q, \mathcal{G}')$ ← customer-meeting-pointAssign($R, \mathcal{G}$)
3. set up e-DARP($q, \mathcal{G}'$) instance
4. $s_{init}$ ← generateInitSol(e-DARP)
5. set the current solution $s=s_{best}=s_{init}$, $i_{imp}=0$, $T=T_{max}$.
6. for $iter = 1 : iter_{max}$
7.     $i_{imp} \leftarrow i_{imp} + 1; s \leftarrow s'$
8.     if $count_{stagnant} = n_{stagnant}$
9.         return $s_{best}$
10.    end if
11.    ls ← $rand(LS)$
12.    $s' \leftarrow$ ls($s$)
13.    if $c(s') < c(s) + T$ and noChargingConflict($s'$) = true
14.        $s'' \leftarrow$ vehicleExchange($s'$)
15.        $s \leftarrow s''$
16.        If $c(s) < c(s_{best})$ and $n_k(s) < |K| + 1$
17.            $s_{best} \leftarrow s; i_{imp} = 0$
18.        end if
19.    end if
20.    if $i_{imp} > 0$
21.        $T \leftarrow T - T_{max}/T_{red}$
22.        if $T < 0$
23.            $T \leftarrow random(0,1) * T_{max}$
24.            if $i_{imp} > n_{imp} * n_k(s_{best})$
25.                $s \leftarrow s_{best}; i_{imp} = 0$
26.            end if
27.        end if
28.    end if
29. end for
30. return $s^*_{best} \leftarrow$ postOptimization($s_{best}$)
```

4.3. Assignment of customers to meeting points

Given customers' origins and the maximum walking distance constraint, we need to determine which meeting points customers should be assigned to considering the trade-off between bus routing costs and customer inconvenience (walking time). The considered customer-to-meeting-point assignment problem is formulated as an MILP of a variant of the capacitated facility location problem. Let $\ell$ denote a layer, and $\mathcal{L}$ denote the set of layers of a layered graph. $G'_\ell$ is the subset of dummy meeting point nodes of layer $\ell$. The objective function of equation (36) minimizes the weighted sum of total customer walking time ($t_{rj}$) and bus travel time ($t_{ij}$) between activated (with positive assigned customers) meeting points. $\lambda_1$ and $\lambda_2$ are the weights in the original objective function (Eq. (1)). $\{\rho_\ell\}$ is non-negative coefficients to be tuned to trade-off between customer walking time and vehicle travel time (see Section 5.2 for the tuning method). $y_{rj}$ and $\theta_j$ are binary variables indicating whether customer $r$ is assigned to meeting point $j$ and whether meeting point $j$ has positive assigned customers, respectively. $z^\ell_{ij}$ (eq. (36)) is an indicator if arc $(i,j)$ on layer $\ell$ is used.

$$\text{Min } \lambda_2 \sum_{r \in R} \sum_{j \in G'_{\ell(r)}} t_{rj} y_{rj} + \lambda_1 \rho_\ell \sum_{\ell \in \mathcal{L}} \sum_{i \in G'_\ell} \sum_{j \in G'_\ell} t_{ij} z^\ell_{ij} \quad (36)$$

Subject to
$$c_{rj} y_{rj} \leq w_{max}, \forall r \in R, j \in G'_{\ell(r)} \quad (37)$$

$$\sum_{j \in G'_{\ell(r)}} y_{rj} = 1, \forall r \in R \quad (38)$$

$$\sum_{r \in R} y_{rj} \leq Q_{max}, \forall j \in G'_{\ell(r)} \quad (39)$$

$$\sum_{r \in R} y_{rj} \leq M\theta_j, \forall j \in G'_{\ell(r)} \quad (40)$$

$$z_{ij}^{\ell} \leq \theta_i, \forall i \in G'_{\ell}, j \in G'_{\ell}, \ell \in \mathcal{L} \quad (41)$$

$$z_{ij}^{\ell} \leq \theta_j, \forall i \in G'_{\ell}, j \in G'_{\ell}, \ell \in \mathcal{L} \quad (42)$$

$$z_{ij}^{\ell} \geq \theta_i + \theta_j - 1, \forall i \in G'_{\ell}, j \in G'_{\ell}, \ell \in \mathcal{L} \quad (43)$$

$$z_{ij}^{\ell} \in \{0,1\}, \forall i \in G'_{\ell}, j \in G'_{\ell}, \ell \in \mathcal{L} \quad (44)$$

$$\theta_j \in \{0,1\}, \forall j \in G'_{\ell(r)}, \forall r \in R \quad (45)$$

$$y_{rj} \in \{0,1\}, \forall j \in G'_{\ell(r)}, \forall r \in R \quad (46)$$

Equations (37)-(38) ensure that each customer is connected to exactly one meeting point within the maximum walking distance. Equation (39) states that the total number of customers assigned to a meeting point cannot exceed the maximum vehicle capacity. Equation (40) ensures the consistency constraint between $y_{rj}$ and $\theta_j$. $Q_{max}$ is the maximum of vehicle capacity. Equations (41)-(43) state that $z_{ij}^{\ell}$ is equal to 1 when both $\theta_i$ and $\theta_j$ are 1 with $i,j \in G'_{\ell}$. $M$ is the number of customers $|R|$. The above MILP problem can be easily solved using a commercial MILP solver to provide a good starting point of the customer-to-meeting-point assignment for the subsequent electric bus routing problem with capacitated charging stations. Note that, to reduce the problem size and solve it efficiently, $G'_{\ell}$ contains only active (i.e. within the maximum walking distance of customers' origins) dummy meeting points. Non-active nodes in $G'$ are trimmed off.

### 4.4. Bus route optimization with charging synchronization constraints

#### 4.4.1. Generation of initial solutions and charging scheduling of the vehicles

We randomly generate $n$ solutions based on the greedy insertion operator. This operator inserts one customer at a time at the least cost and feasible position of that route by checking the time window, maximum ride time, precedence, and vehicle capacity constraints using the eight-step evaluation scheme (Cordeau and Laporte, 2003). Afterwards, the energy constraints resulting from that insertion are checked. If violated, charging operations (insertion of visits to chargers) are scheduled. The charging scheduling algorithm is described in Algorithm 2. The charging scheduling of the current route first identifies energy feasible positions after which the vehicle can go to recharge (line 2 in Algorithm 2). These occur after leaving the depot, or after leaving a transit station, because we assume that charging operations are only allowed when there are no passengers onboard. Note that direct connection from one charger to another is forbidden. Given a set of feasible charging positions for the current route, we calculate the forward slack time (line 5 in Algorithm 2). This slack time allows to calculate the time available for the vehicle's recharging without violation of the time windows of the remaining route. To reduce the risk of conflicts (charging duration overlap) with charging operations of other vehicles (charging synchronization), we adopt a mixed randomization strategy. First, the charging position is randomly selected among the list of candidate charging positions of the current route (line 3 in Algorithm 2). Given the selected insertion position, a charger is selected with a greedy strategy (i.e. select a charger over all chargers with the least charging operation time, including the access, egress and charging times). Given the slack time for recharging and the desired charging duration, a feasible starting time for recharging is randomly selected within a reduced interval without compromising the desired charging duration (lines 7-10 in Algorithm 2). If the vehicle cannot be recharged to the desired level for the current charger insertion, a subsequent recharging visit is inserted at the next feasible location in the current route (lines 12-14 in Algorithm 2). When the above charging scheduling attempt fails, a charging recovery procedure is invoked by restarting the charging insertion from the first feasible position (lines 19-22 in Algorithm 2). If no feasible charging insertion can be found, the route is considered infeasible. The output is a sequence of charging operations, each containing the insertion position, charger node to be visited, and the start and termination times of recharge. The charging conflict check with other vehicles is invoked when the current temporary solution is promising (line 13 in Algorithm 1). The charging station occupancy state of a solution is handled using a discretization scheme. This is implemented as a binary matrix where each line represents the occupancy of a charger over the planning horizon with a discretized time interval of 10 seconds. By recording the charging schedules of the solution routes on this matrix, this conflict check can be conducted in O(1) time.

**Algorithm 2. Vehicle charging scheduling algorithm.**

1. **Input:** an energy-infeasible temporary route $r = \{v_0, \ldots, v_n\}$.
2. Set $success = false, repair = false$. Compute the cumulative energy consumption from the depot to each node of $r$. Identify the first infeasible node $v^*$ where the state of charge is insufficient to reach that node (i.e. inferior to the $E_{min}^k$).
3. Let $loc_r = \{(v_{i-1}, v_i)\}_{i=1,\ldots,p}$ be a list of candidate positions for inserting a recharge operation where a possible recharge visit to a charger location is between $v_{i-1}, v_i$. Randomly select a position between 1 and $p$ to insert the first recharging visit.
4. Choose a charger among the list of all chargers with the least charging operation time, including the access, egress and charging times to connect $v_{i-1}$ and $v_i$.
   Let $v_s$ denote this charger node. Compute the amount energy of recharge $\Delta_{E_r} = E_{min}^k - E_{init}^k - ec_r$, where $ec_r$ is the total energy consumption of $r$ including the recharge visit. Compute the charging times of $\Delta_{E_r}$ given the power of the charger to be visited, denoted as $\delta_{\Delta_{E_r}}$.
5. Compute the allowed delays of the service start time at $v_{i-1}$ as $\delta_{i-1}^{delay} = min\{F_{i-1}, \sum_{i-1<p<n} W_p\}$ for all candidate positions, where $F_{i-1}$ and $W_{i-1}$ are the forward time slack and the waiting time at node $v_{i-1}$, respectively (Cordeau and Laporte, 2003).
6. Try inserting $v_s$ after $v_{i-1}$ as follows. Compute $\delta_{access} = T_{v_{i-1}, v_s} + T_{v_s, v_i} - T_{v_{i-1}, v_i}$.
7. **If** $\delta_{\Delta_{E_r}} + \delta_{access} \leq \delta_{i-1}^{delay}$
8.     Compute the recharge start time $B_{v_s} := D_{v_{i-1}} + T_{v_{i-1}, v_s} + rand(\delta_i^{delay} - \delta_{\Delta_{E_r}})$, where $D_{v_{i-1}}$ is the
9.     departure time, and $rand()$ generates a random real number between its input range. Compute the
10.     terminate time of recharge $T_{v_s}^{end} = B_{v_s} + \delta_{\Delta_{E_r}}$. Set $success = true$.
11. **else**
12.     Compute $B_{v_s} := D_{v_{i-1}} + T_{v_{i-1}, v_s}$ and $T_{v_s}^{end} := B_{v_s} + (\delta_{i-1}^{delay} - \delta_{access})$. Update the remaining amount
13.     energy to be recharged after $v_s$. Perform the same charging scheduling procedure for the remaining
14.     candidate positions. If successful, set $success = true$.
15. **end**
16. **If** $success = true$
17.     return the charging schedules of $r$, i.e. $\{v_{i-1}, v_s, B_{v_s}, T_{v_s}^{end}\}$.
18. **else**
19.     If $repair = false$
20.         Set $repair = true$. Restart the charging scheduling from the first candidate position, i.e. $v_{i-1}$ = the
21.         depot, go to line 6.
22.     end
23. **end**

### 4.4.2. Local search operators

The local search operators need to be designed by considering their complementarity in diversifying the searched neighborhood from the current solution (Arnold and Sörensen, 2019) and the specificity of the problem at hand. As we allow customer requests to be rejected, unserved customers are managed in a pool, which is regarded as a **virtual route**, allowing customers [2] to be removed from the vehicles. Note that charging schedule updating is applied at the end of each local search operation. Given the layered graph structure, we can efficiently screen out infeasible insertion positions by checking whether the layer of a customer to be inserted is (in)compatible with the layer of the current inserted position of the route. This conflict check can be done in O(1) based on a lookup table of compatibility information for the layered graph. However, the computational time savings remain marginal. The effectiveness of the layered graph structure is mainly demonstrated in solving the customer-to-meeting-point assignment model. The computational results for the effectiveness of using the layered graph structure are shown in Appendix C. We propose eight local search operators as follows.

---
[2] Note that a customer in the bus route optimization procedure denotes the requests assigned to the same dummy meeting point based on the customer-to-meeting-point of Eqs. (36)-(46).

- **Relocate ensemble**: To avoid repetitive local search operations applied to the same neighborhood, we randomly select a relocation operator from the following two relocation operators: greedy relocation and worst relocation. The two relocation operators search different parts of a solution neighborhood. Greedy relocation randomly removes a customer from their current route and reinserts the customer to the least cost position of the current solution. Worst relocation removes the worst customer (i.e., the most expensive with respect to the objective function) of a randomly selected route and reinserts the customer into the least cost position of the current solution.
- **Destroy-repair**: To shake the current solution significantly, we adopt the idea of destroy-repair idea used in the Adaptive Large Neighborhood Search heuristic (ALNS) (Ropke and Pisinger, 2016). Different from ANLS, our destroy-repair operator removes $n^{remove}$ based on a randomly selected removal operator and reinserts the unserved customers using a randomly selected repair operator. The number $n^{remove}$ is selected randomly as $1 \leq n^{remove} \leq \min(n^{max}, \delta|R|)$, where $n^{max}$ is a pre-defined maximum number of customers to remove. $\delta$ is a coefficient between 0.2 and 0.5. Based on our preliminary analysis, we set $n^{max} = 60$, and $\delta = 0.275$ (Lutz, 2014). Five removal operators and two repair operators are used as follows.

*Destroy operators*
**Random removal**: Randomly remove $n^{remove}$ customers from the candidate list (i.e. all customers on the current solution $s$).
**Worst removal**: Let $c(s)$ and $c(s_{-i})$ denote the cost with and without customer $i$ on the current solution $s$ Compute the cost of the customer $i$ as $\bar{c}(s_{-i}) = c(s) - c(s_{-i})$ for all $i \in s$. Sort all customers in $s$ in descending order of $\bar{c}(s_{-i})$. To avoid repeatedly selecting the same worst customers to remove, a randomness strategy is used where the degree of randomness is controlled by a parameter $p$ (Ropke and Pisinger, 2016). Compute the determination parameter $y^p$ where $y$ is a random number between 0 and 1, and $p$ is the degree of randomness parameter. Select $j$-th customer to remove on the sorted list, where $j = \max(1, \lfloor y^p n^{remove} \rfloor)$. Continue until $n^{remove}$ customers are removed. The reader is referred to Ropke and Pisinger (2016) for more detailed description.
For the following related-removal operators, we describe the **removal algorithm** as follows. First, compute the relatedness indicators (described below) for all pairs of customers. We remove $n^{remove}$ customers from $s$ as follows. First, a customer is removed randomly from $s$ and added to the unserved pool. On the subsequent iterations, select initially a customer $i$ randomly from the unserved pool, then sort the relatedness to customer $i$ for all customers in $s$ in ascending order. Compute the determination parameter $y^p$ as above. Select the $j$-th customer to remove from the sorted list, where $j = \max(1, \lfloor y^p n^{remove} \rfloor)$. Continue until $n^{remove} - 1$ customers are removed from $s$. We set $(p^{worst}, p^{dist}, p^{tw}, p^{Shaw},)$ as (3, 6, 6, 6) (Ropke and Pisinger, 2016).
**Distance-related removal**: Compute the distance relatedness for all pairs of customers as follows (adapted from Lutz, 2014).

$$r_{dist}(i,j) = dist(i,j) + dist(n+i, n+j) \tag{47}$$

where $dist(i,j)$ is the normalized travel time (i.e. travel time divided by a maximum travel time value) between the pickup locations of customers $i$ and $j$. While $dist(n+i, n+j)$ corresponds to the normalized travel time of their drop-off locations. Then remove $n^{remove}$ customers based on the removal algorithm.
**Time-window-related removal**: Compute time-window relatedness for all pairs of customers as follows.

$$r_{tw}(i,j) = |\tilde{l}_i - \tilde{l}_j| + |\tilde{l}_{n+i} - \tilde{l}_{n+j}| \tag{48}$$

where $\tilde{l}_i$ is the normalized latest starting times of service at vertex $i$ (i.e. $l_i/T$ where $T$ is the planning horizon). Then remove $n^{remove}$ customers based on the removal algorithm.
**Shaw removal**: Compute the Shaw relatedness (adapted from Ropke and Pisinger (2016)) for all pairs of customers as follows.

$$r_{Shaw}(i,j) = \varphi r_{dist}(r_i, r_j) + \chi r_{dist}(r_i, r_j) + \psi|\tilde{g}_i - \tilde{g}_j| \qquad (49)$$

where $\tilde{g}_i$ is the normalized capacity demand of customer $i$ (with respect to max{ $g_i$ } for all $i \in G'$ where $g_i$ is the number of passengers associated with $i$. Then remove $n^{remove}$ customers based on the removal algorithm. We set ($\varphi, \chi, \psi$) as (9, 3, 2) (Ropke and Pisinger, 2016).

*Repair operators*
**Greedy insert**: Insert the customers of the unserved pool one by one to the least cost positions of the current solution.
**Regret insert**: Insert the customers of the unserved pool one by one based on the regret heuristic of Ropke and Pisinger (2016). Note that we implement Regret-2 and Regret-3 heuristics, which are selected randomly to achieve more diverse insertions.

- **Two-opt***: Two routes are randomly selected. Identify the candidate arcs to be removed from each route, i.e. the load of the vehicle is zero on these arcs (Parragh et al., 2010). Remove one candidate arc from each route and recombine the first part of the first route with the remaining part of the second route, and vice versa. The feasible one with the best cost improvement is retained.
- **Two-opt**: Reverses the order of visiting the nodes of a segment of the current route, in a sequential manner along the current route. First, the length of the segment is randomly determined between 2 and 4. A node is randomly selected and the end node of the segment can be identified. The feasible one with the best cost improvement is kept.
- **Exchange-segment**: Randomly select two routes. Identify candidate arcs with zero passenger load over the entire current route. The segments between two consecutive candidate arcs of a route are the swappable segments. Exchange two such segments, sequentially along the route, until a feasible and improved route is found (i.e., the vehicle travel time savings after the exchange is greater than the current threshold T).
- **Exchange-customer**: Swap two customers on two randomly selected routes. The pickup (drop-off) position of the customer on the first route needs to be reinserted at the pickup (drop-off) position of the customer on the second route (Braekers et al., 2014). If successful, the removed customer of the second route is reinserted into the first route at any feasible position. If unsuccessful, the customer is reinserted on another randomly selected route until a feasible insertion is found. This operator is applied sequentially along the routes until an improved and feasible exchange is found.
- **Four-opt**: Remove four successive arcs (three successive nodes) from the current route and find a feasible and improved one among all possible combinations of the removed segment in a sequential search along the route. The feasible and improved one with the least cost is retained.
- **Create**: if the number of used vehicles is smaller than the fleet size and the pool of unserved customers is not empty, create a new route by randomly inserting an unserved customer on an unused vehicle.

Note that at the end of each local search operator, we update the vehicle charging schedule by applying the vehicle charging scheduling procedure (Algorithm 2). If the resulting routes satisfy all the constraints, the updated solution s' is kept (line 12 in Algorithm 1).

4.5. Post-optimization

If the best solution obtained by the DA algorithm contains unserved customers, the customers of the layers containing unserved customers are reassigned and the subroutes on the affected layers are reoptimized to satisfy the time window, maximum ride time, and vehicle capacity constraints. The reassignment problem is solved by a matheuristic that formulates a conventional MILP problem for reassigning customers on the same layer that contains unserved customers (Appendix B). This matheuristic aims to find improved partial routes for each of such layers that minimize the total vehicle routing time and the total penalty of unserved customers without violating constraints. A short computation time limit is applied to find an optimal solution. If the solution reduces the initial number of unserved customers on the current layer, the corresponding subroutes are updated. Finally, the charging schedules of the modified routes are updated and checked for charging conflicts. If the new solution is feasible and has an improved cost, the best solution is updated. The details of the matheuristic are described in Algorithm 3 (see Appendix B).

## 5. Computational study

In this section, we present computational results for the proposed algorithm and compare its performance with exact solutions obtained by a state-of-the-art MILP solver. First, we present the test instances and tune the algorithmic parameters. Then, the performance of the hybrid metaheuristic is evaluated. A sensitivity analysis is performed to investigate the impact of different system parameters. To simulate the application to real-world cases, a case study mimicking a feeder service to serve a central city train station is studied. The sensitivity analysis to determine system parameters is analyzed.

5.1 Test instance generation

We implement a scenario generator that can create different test instances based on several system parameters: number of customers, number of transit stations, number of chargers, meeting point separation distance, maximum walking distance, number of punctuated feeder services per transit station, distance between customer locations, and demand distribution through the planning period: parametrized from uniform distribution to peaked normal distribution. We consider two scenarios corresponding to peak (P) and off-peak (OP) demand profiles. Scenario P simulates a peak-hour situation where customers' desired arrival times are concentrated around a peak hour, while OP reflects the opposite situation when customers' desired arrival times are uniformly spread over a longer operating period. Each test instance contains bus system supply and network data (locations of the depot, meeting points, and transit stations), and charging infrastructure, and randomly generated customer demand information (origins, destinations (drop-off transit stations), desired arrival time at destinations). To get a nuanced performance of the algorithm, we generate 10 instances in each scenario, ranging from 10 to 100 customers. These test instances have a single vehicle depot, two train stations, and four chargers as shown in Figure 3. Meeting points are generated as a grid with a separation distance of 1 km and customers maximum walking distance is 1.5 km. Punctuated services are provided for the two train stations with three services per hour starting from 6:00 (6:10) and ending at 10:00 (9:50), respectively. In total, there are 26 layers with 25 to 49 activated meeting points per layer (meeting points within the maximum walking distance of the customers). Each customer may have up to 7 meeting points within their walking distance (Figure 3). Consequently, the possible customer-to-meeting-point assignment combinations are very large, providing non-trivial experiments to test the performance of the algorithm. We consider two types of vehicles with different passenger capacity, battery capacity, and energy consumption rate. The instance name cxx means that there are xx customers in that instance. In order to force scenarios where vehicles need to recharge, initial battery levels of vehicles are set as low as 20%, 30%,…, and 80% of the battery capacity. The number of vehicles is determined by considering approximately 70% of vehicle occupancy. The average initial battery state for both scenarios is 26.4%. We assume that, when returning at the depot, the minimum battery level is 10% of its capacity. An overview of the experiment is shown in Table 1.

**Table 1. Overview of the experiment and algorithmic parameter settings**

| Parameter | Value | Parameter | Value |
|---|---:|---|---:|
| Number of type of vehicles (buses) | 2 | Charging power | 0.83kWh/minute |
| Number of meeting points within customer's walking distance | 25 to 42 | $E_{init}^k$ | 20%, 30%,…,80% of $\bar{B}_k$ |
| Number of punctuated service per transit station | 12 or 13 | Energy consumption rate of vehicles | 0.24 kWh/km and 0.29 kWh/km |
| Number of transit stations | 2 | Number of customers | 10,20,…,100 |
| Number of chargers | 4 | Walking speed | 0.085 km/minute |
| Number of vehicles | 2 to 6 | Detour factor | 1.5 |
| Passenger capacity of the vehicles | 10 or 20 | $\lambda$ | 1 |

| Vehicle speed | 0.83 km/minute | $\omega$ (penalty of one unserved customer)* | 40 |
| --- | --- | --- | --- |
| Battery capacity of vehicles ($\bar{B}_k$) | 35.775 kWh and 53.70 kWh | Maximum walking distance of customers | 1.5 km |
| $E_{max}^k (E_{min}^k)$ | $0.8\bar{B}_k$ ($0.1\bar{B}_k$) | $u_i$ (service time) | 0.5 minutes |

*based on our preliminary experiments. The characteristics of electric buses is adapted from Volkswagen's 8-seat 100% electric Tribus with 35.8 kWh (https://www.tribus-group.com/zero-emission-volkswagen-e-crafter-electric-wheelchair-minibus/)

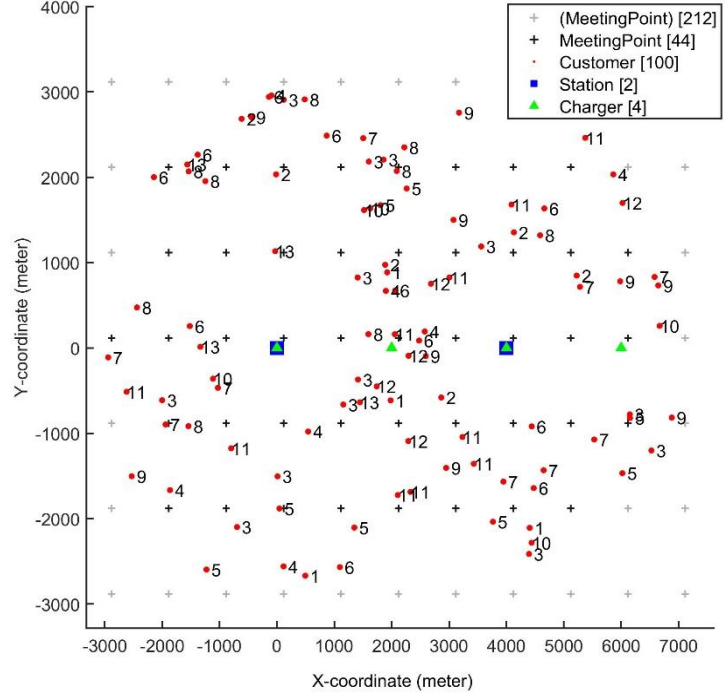

**Figure 3. Example of a test instance with 100 randomly generated customer locations. Number next to each customer shows desired transit departure. The greyed out (Meeting point) are beyond the maximum walking distance of any customer.**

5.2. Algorithmic parameter settings

The parameters used in Algorithm 1 include the customer-to-meeting-point assignment weight ($\rho$) and several parameters used for the DA algorithm. We first present the sensitivity analysis of the parameters of the DA algorithm. Following Braekers et al. (2014), the parameters of the DA algorithm and the associated discrete values to be tested are listed as follows.

- $t_{max}$: A user-defined coefficient to determine the maximum threshold value $T_{max}$ to accept worsen temporary solutions, i.e. $T_{max} = t_{max} T_{allbus}$ where $T_{allbus}$ is the average bus travel time between all pick-up and drop-off nodes in the layered graph of a problem instance. 8 values of $t_{max}$ are tested, i.e. $t_{max} \in \{0.3, 0.6, 0.9, 1.2, 1.5, 1.8\}$.
- $T_{red}$ : Threshold reduction factor for reducing the threshold value ($T \coloneqq T - \frac{T_{max}}{T_{red}}$) (line 20 in Algorithm 1). Eight values of $T_{red}$ are tested, i.e. $T_{red} \in \{100, 200, 300, 400, 500, 1000\}$.
- $n_{imp}$: Restart parameter, set as $n_{imp} \in \{100, 200, 300, 400, 500, 600\}$.
- $Iter_{max}$: Maximum number of iterations, set as $Iter_{max} \in \{25, 50, 100, 150, 200, 300, 400\}$.
- $n_{stagnant}$: Maximum stagnation multiplier, set as $n_{stagnant} \in \{25, 50, 100, 150, 200, 250, 300\}$

We generate 10 random test instances with random numbers of requests between 20 and 100 (see Table 3). We solve each instance 5 times with a constant $\rho$ to get the average performance for each tested parameter setting. It takes a few days to complete the experiments on a high-performance machine.

To evaluate the performance of a solution, a gap $(Z - Z^*)/Z^*$ is calculated with respect to the best solution found for each instance. Based on our preliminary analysis, $t_{max}$, $T_{red}$ and $n_{imp}$ are initialized as 0.6, 500 and 200, respectively. We first tune $Iter_{max}$ without considering $n_{stagnant}$. Given the tuned $Iter_{max}$, we vary $n_{stagnant}$ to analyze its sensitivity and a best value. Since the effects of $t_{max}$, $T_{red}$, and $n_{imp}$ on the performance of the algorithm are not independent (Braekers et al., 2014), the best parameter setting needs to consider all possible combinations of the tested values of these parameters. The left side of Figure 4 shows the effect of $Iter_{max}$ on the average gaps over the 10 test instances, all other parameters being equal. The best $Iter_{max}$ is found at 300k iterations with an average gap of 0.06%. However, the computation time is proportional to $Iter_{max}$. For the remaining experiments, we decide to use $Iter_{max} = 100k$ which is a good trade-off between computation time and solution quality (gap=0.17%). On the right side of Figure 4, the best $n_{stagnant}$ is found at 200 with the average gap of less than 0.15%. Increasing $n_{stagnant}$ further to 300 significantly increases the computation time with little improvement in the gap value. The sensitivity analysis of $t_{max}$, $T_{red}$, and $n_{imp}$ is shown in Table 2. We can observe that in general the gap increases with $t_{max}$, $T_{red}$, and $n_{imp}$. From the result of 512 combinations, the minimum gap 0.13% is found with $t_{max} = 2.1$, $T_{red} = 200$ and $n_{imp} = 100$. Therefore, we will use these values with $n_{iter} = 100$k and $n_{stagnant} = 200$ for the computational studies in the following sections.

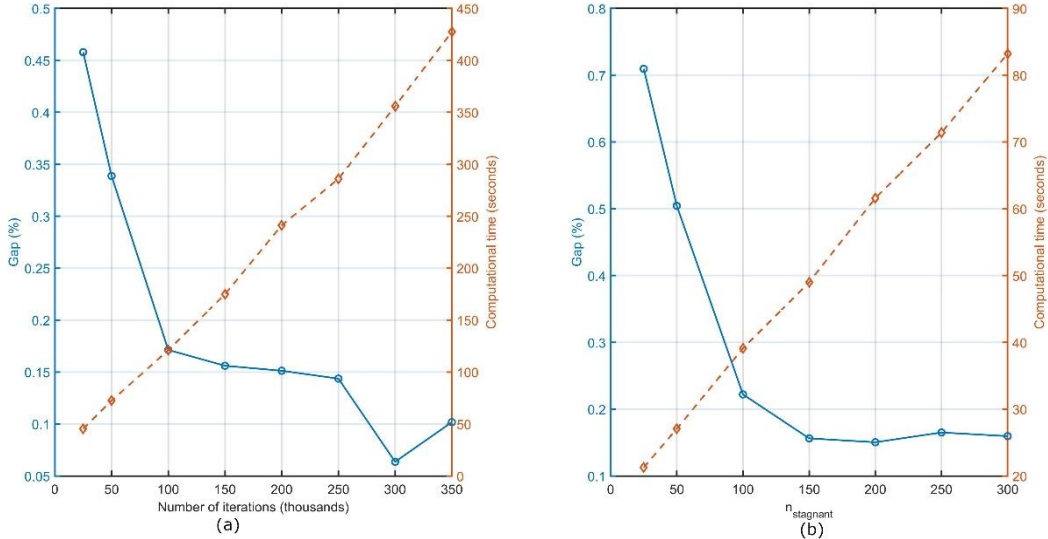

**Figure 4.** Impact of $Iter_{max}$ (a) and $n_{stagnant}$ (b) on the performance of the hybrid metaheuristic.

**Table 2.** Sensitivity analysis on $t_{max}$, $T_{red}$ and $n_{imp}$.

| $t_{max}$ | 0.3 | 0.6 | 0.9 | 1.2 | 1.5 | 1.8 | 2.1 | 2.4 |
|---|---|---|---|---|---|---|---|---|
| Avg. gap | 0.54% | 0.61% | 0.60% | 0.64% | 0.62% | 0.67% | 0.69% | 0.67% |
| $T_{red}$ | 100 | 200 | 300 | 400 | 500 | 1000 | 2000 | 3000 |
| Avg. gap | 0.58% | 0.61% | 0.62% | 0.61% | 0.61% | 0.63% | 0.68% | 0.69% |
| $n_{imp}$ | 100 | 200 | 300 | 400 | 500 | 600 | 1000 | 2000 |
| Avg. gap. | 0.36% | 0.49% | 0.56% | 0.65% | 0.66% | 0.68% | 0.79% | 0.84% |

Regarding the impact of $\boldsymbol{\rho} = \{\rho_\ell\}_{\ell \in \mathcal{L}}$, it depends on the characteristics of each test case. Different factors such as demand intensity, spatial distribution of requests, maximum walking distance of customers, number of vehicles, capacity of vehicles, location of meeting points and transit stations may affect the best value of $\boldsymbol{\rho}$ for each test instance. Figure 5 shows the effect of $\boldsymbol{\rho}$ on the performance of the algorithm on three test instances with 28, 52, and 88 customers. It can be observed that the relationship between $\rho$ and the algorithm performance is not unimodal. Thus, the classical golden section search method is not relevant to find the best value of ρ. The best **ρ** tends to be smaller (0.2) for larger instances to minimize the customer walking time in Eq. (36). We conduct a systematic $\rho$ search experiment over the above 10 test instances using a 2-step approach as follows. First, set 10 test values of $\rho$ (i.e. $\rho_\ell = \rho$ for all $\ell \in \mathcal{L}$) as $\rho \in P = \{0.2, 0.4, \ldots, 2.0\}$. Get the preliminary best $\tilde{\rho}$ (best solution found) over P for each test instance. Then run two tests to get the best $\rho^*$ on the neighborhood of $\tilde{\rho}$, i.e. $\tilde{\rho} \pm 0.1$. For instances with unserved users, we tune the layer-specific $\rho_\ell$ with 5 additional runs on a second step. Initially set $\rho^0 = \rho^*$ for all $\ell \in \mathcal{L}$. For the subsequent iterations, increase the value of $\rho_\ell$ for $\ell \in \bar{\mathcal{L}}$ with $\delta$ (e.g. 1.5), where $\bar{\mathcal{L}}$ is the subset of layers with non-empty unserved users. We retain $\{\rho_\ell^*\}_{\ell \in \mathcal{L}}$ of the best of the 5 runs. It can be observed that for instances with more than 50 customers, the best $\rho^*$ is no larger than 0.5, while for small instances, the best $\rho^*$ can be larger than 1. This suggests that a narrowed search range between 0 and 0.5 is sufficient to find good $\rho^*$ for the cases with more than 50 requests. In practice, it is not a problem to find the best $\rho^*$ with more test values for small instances since the computation time is in the order of seconds. On the other hand, for larger instances, a budget of 5 runs can be allocated to find the best solution using 5 test values of $\rho$, i.e. $\rho \in \{0.2, 0.4, 0.6\}$ with two additional test values on the neighborhood ($\pm 0.1$) of $\tilde{\rho}$. If there are unserved users, additional 3 runs are necessary to adjust $\rho$ for layers with unserved users.

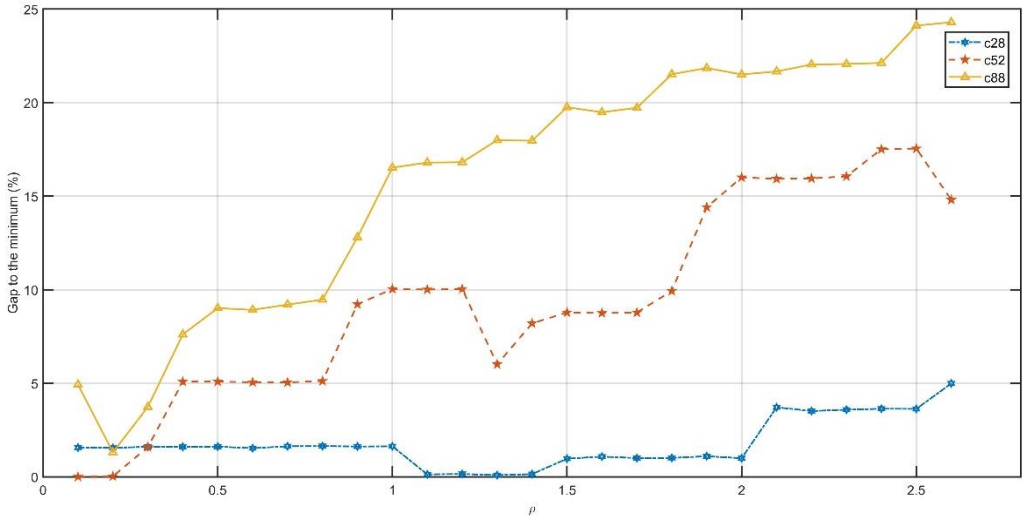

**Figure 5. Impact of $\rho$ on the performance of the hybrid metaheuristic algorithm.**

**Table 3. Best values of $\rho$ on the test instances.**

| Instance | c24 | c28 | c30 | c34 | c40 | c48 | c52 | c72 | c82 | c88 |
|---|---|---|---|---|---|---|---|---|---|---|
| Number of vehicles | 2 | 2 | 3 | 3 | 3 | 3 | 3 | 5 | 5 | 6 |
| Average $E_{init}^k$ (% of $\bar{B}_k$) | 20 | 20 | 23 | 23 | 23 | 23 | 23 | 28 | 28 | 32 |
| Best $\rho$ | 0.7 | 1.2 | 0.7 | 0.1 | 0.6 | 1.1 | 0.2 | 0.4 | 0.5 | 0.1 |

Remark: $\rho$ is the same for each layer (Eq. (36)).

### 5.3. Computational results

We test the performance of the algorithm on the 20 test instances in peak and off-peak scenarios using the tuned parameters for the DA algorithm. The algorithm performance is compared with the solution obtained by a state-of-the-art MILP solver (Gurobi, version 10) with a computation time limit of 4 hours. Our algorithm and the MILP model are both implemented using the Julia programming language. We run the experiments on a laptop with an 11th generation Intel(R) Core(TM) i7-11800H, 16 logical processors and 64 GB of memory. It is worth noting that Gurobi v10 is a parallel MIP solver that uses up to 16 threads by default, while the proposed algorithm runs on a single thread. It is quite obvious that by implementing our algorithm using parallel computing technique, the computation time could be significantly reduced. The MILP solutions obtained by Gurobi are shown in Table 4. For each instance, the MILP results give the best feasible solution found within 4 hours, along with the lower bound. The third column shows the number of unserved customers. Note that each unserved customer is associated with a high penalty (40/customer). This has a large impact on the value of the objective function. The solver can obtain (near) optimal solutions for small instances of 10 customers. Within 4 hours of computation time, the number of unserved customers is zero for 14/20 test instances, with one unserved customer for 6 of the P scenarios. The charging time of the best solution is shown next to the number of unserved customers. We observe that charging operations are present on all test instances.

  For the hybrid metaheuristic, the results are based on the average of 5 runs with random seeds to account for the random elements within the algorithm. For each instance, we report the best objective function value (Best obj.), the gap between the average objective function value and that of the best solution to the best known solution (BKS) found by the solver. The last four columns report the number of unserved customers, charging time of the best solution, the average computation time (per run). Overall, the average and best gaps for the 20 test instances are 1.80% and 1.69%, respectively. For the OP scenario, the algorithm gets the same solution (with 10 customers) and small gaps for the other instances. Overall, the average and best gaps are 0.67% and 0.49%, respectively. We observe the average gaps to the BKS are within the range of -2.82% and 3.69%. We observe that higher gaps are found for c50op with the worst average gap of 3.69%. Similar results can be found for scenario P with a slightly higher average and best gap of 2.93% and 2.89%, respectively. The overall average computation time of the algorithm is less than 1 minute per run. For OP scenario, all customers are served, while for the P scenarios only c30p and c50p result in an unserved customer (an improvement over the 6 MILP scenarios with unserved customers). For the largest instance with 90 or 100 requests, the OP scenario requires more computation time (~2.5 minutes) due to its larger problem size (the number of active layers is much higher compared to the P scenario for the same problem size). The results show the efficiency and good solution quality obtained by the proposed algorithm. Note that the worst gap is found for the c50p case with an average gap of 11.19%; this also appears to be a difficult instance for the MILP in terms of the BKS and gap to the lower bound. These results indicate that for some difficult instances, more computational effort are needed to explore potential customer-to-meeting-point assignment alternatives.

**Table 4. Computational results obtained using the Gurobi solver and the hybrid metaheuristic algorithm on the test instances.**

| Instance | MILP | | | | Hybrid metaheuristic | | | | | |
|---|---|---|---|---|---|---|---|---|---|---|
| | Best known solution (BKS) | Gap to the lower bound | UC* | CT* (min) | Best obj. | Gap** (avg.) | Gap** (best) | UC* | CT* (min) | CPU (s) |
| c10op | 107.91 | 4.29% | 0 | 4.9 | 107.91 | 0.00% | 0.00% | 0 | 4.9 | 8 |
| c20op | 232.57 | 15.14% | 0 | 17.9 | 234.15 | 0.69% | 0.68% | 0 | 18.2 | 12 |
| c30op | 327.76 | 13.90% | 0 | 15.4 | 334.26 | 2.10% | 1.98% | 0 | 16.5 | 28 |
| c40op | 413.96 | 23.49% | 0 | 29.2 | 416.88 | 0.77% | 0.71% | 0 | 31.1 | 26 |
| c50op | 566.74 | 31.74% | 0 | 44.1 | 586.80 | 3.69% | 3.54% | 0 | 42.4 | 30 |
| c60op | 596.35 | 24.05% | 0 | 22.8 | 602.63 | 1.24% | 1.05% | 0 | 26.1 | 84 |
| c70op | 702.64 | 25.09% | 0 | 32.7 | 694.31 | -0.83% | -1.18% | 0 | 33.4 | 76 |
| c80op | 842.29 | 32.13% | 0 | 40.0 | 815.56 | -2.82% | -3.17% | 0 | 43.1 | 85 |
| c90op | 816.47 | 29.05% | 0 | 25.1 | 816.17 | 0.30% | -0.04% | 0 | 31.5 | 149 |
| c100op | 967.48 | 27.79% | 0 | 33.2 | 980.16 | 1.56% | 1.31% | 0 | 27.2 | 142 |
| Average | 557.41 | 22.67% | 0 | 26.5 | 558.88 | 0.67% | 0.49% | 0 | 27.4 | 64 |
| c10p | 112.02 | 16.81% | 0 | 2.8 | 112.55 | 0.47% | 0.47% | 0 | 4.9 | 6 |
| c20p | 283.40 | 40.72% | 1 | 13.5 | 311.26 | 9.84% | 9.83% | 0 | 12.6 | 10 |
| c30p | 341.11 | 33.48% | 1 | 10.7 | 361.15 | 5.88% | 5.88% | 1 | 8.9 | 76 |
| c40p | 421.79 | 39.10% | 1 | 15.7 | 441.80 | 4.74% | 4.74% | 0 | 13.2 | 15 |
| c50p | 606.48 | 46.82% | 1 | 31.9 | 674.19 | 11.19% | 11.16% | 1 | 27.3 | 66 |
| c60p | 584.07 | 34.47% | 0 | 11.6 | 628.98 | 7.86% | 7.69% | 0 | 3.6 | 47 |
| c70p | 684.36 | 34.71% | 1 | 7.7 | 689.05 | 0.75% | 0.69% | 0 | 11.4 | 45 |
| c80p | 872.50 | 43.21% | 1 | 20.7 | 831.52 | -4.65% | -4.70% | 0 | 8.6 | 51 |

| | | | | | | | | | |
|---|---|---|---|---|---|---|---|---|---|
| c90p | 811.70 | 38.11% | 0 | 6.9 | 761.15 | -6.22% | -6.23% | 0 | 18.7 | 85 |
| c100p | 943.47 | 34.25% | 0 | 14.0 | 937.87 | -0.53% | -0.59% | 0 | 14.8 | 61 |
| Average | 566.09 | 36.17% | 0.60 | 13.5 | 574.95 | 2.93% | 2.89% | 0.2 | 12.4 | 46 |
| Overall average | 561.75 | 29.42% | 0.30 | 20.04 | 566.92 | 1.80% | 1.69% | 0.1 | 19.9 | 55 |

*UC: Number of unserved customers; CT: charging time
** Gap to the best solutions found by Gurobi at a computation time of 4 hours. Charging time of the obtained solutions is measured in minutes. CPU time is measured in seconds.

The results of Table 4 use the post-optimization procedure. To explicitly test its performance, in Table 5 we compare the computational results with and without applying post-optimization. We only report the results where there are unserved customers without applying post-optimization to verify the effectiveness of this post-optimization (Table 5). Note that in Table 4 the best objective values are reported, while in Table 5 the average objective values are reported. The results are based on an average of 5 runs. Post-optimization is shown to both improve the solution quality (2/3) and reduce the number of unserved customers to half or zero. The average gap is significantly improved from 12.96% to 8.97%. In terms of computation time, the average time for the algorithm with post-optimization is significantly higher on average (+40 seconds on average), depending on the problem size of the MILP formulation for the post-optimization.

**Table 5. Computational results for the post-optimization procedure of the hybrid metaheuristic algorithm.**

| | Algorithm **without** the post-optimization | | | | Algorithm **with** the post-optimization | | | |
|---|---|---|---|---|---|---|---|---|
| Instance | Avg. obj. Value | Gap* | Num. of unserved customers | CPU (s) | Avg. obj. value | Gap* | Num. of unserved customers | CPU (s) |
| c20p | 335.03 | 18.22% | 3.0 | 6 | 311.29 | 9.84% | 0.0 | 10 |
| c30p | 361.19 | 5.89% | 1.0 | 14 | 361.16 | 5.88% | 1.0 | 76 |
| c50p | 696.09 | 14.78% | 2.0 | 14 | 674.34 | 11.19% | 1.0 | 66 |
| Average** | 464.11 | 12.96% | 2.0 | 11 | 448.93 | 8.97% | 0.7 | 51 |

* Gap to the best solutions found by Gurobi with 4 hours of computation time
**Average over the instances with unserved customers

Note that in practice, vehicles are usually fully charged at the beginning of the day. In this setting, the problem may be easier to solve because the vehicles do not need to be recharged. Table 6 compares the performance of the algorithm with high and low initial battery levels for the 20 instances. We observe that the hybrid metaheuristic outperforms the BKS for the off-peak scenario with high battery levels (easier to be solved compare to the opposite cases). The computation time is less than 17 seconds on average, compared to about 55 seconds on average for the case where $E_{init}^k$ is set to low.

**Table 6. Impact of initial battery level of vehicles on the performance of the hybrid metaheuristic algorithm.**

| Avg. $E_{init}^k$ | MILP | | | | Hybrid metaheuristic | | | |
|---|---|---|---|---|---|---|---|---|
| | Scenario | Avg. BKS | Avg. gap to the lower bound | Charging time | Avg. obj. value | Gap (average) | Gap (best) | CPU (s) |
| Low* | Off-peak | 557.41 | 22.67% | 26.5 | 560.20 | 0.67% | 0.49% | 64 |
| | Peak | 566.09 | 36.17% | 13.5 | 575.22 | 2.93% | 2.89% | 46 |
| High** | Off-peak | 549.21 | 19.34% | 0.0 | 532.23 | -1.70% | -1.81% | 9 |
| | Peak | 573.78 | 35.88% | 0.0 | 562.68 | 1.46% | 1.46% | 24 |

* $E_{init}^k$ is set as 20%, 30%,… with an overall average of 26.4%.

**$E_{init}^k = 80\%$ for all vehicles. BKS is the best solution solved by Gurobi in 4 hours. The result of the hybrid metaheuristic is based on the average of 5 runs. Both sets of experiments use the same algorithmic parameters except $\rho$. Average gap is calculated as the same way as Table 4.

Finally, we compare the stand-alone effectiveness of the metaheuristic for solving the second-stage bus routing problem with charging synchronization and partial recharge problem on the 20 test instances with up to 100 requests. The new experiments utilize the customer-to-meeting-point assignment results as input and then solve the bus routing optimization problem with both Gurobi solver and the metaheuristic (without applying post-optimization for customer reassignment). The MILP formulation of the second stage problem is a heterogeneous electric dial-a-ride problem with charging synchronization and partial recharge. The MILP formulation and computational results are provided in Appendix D. The results show that the proposed metaheuristic outperforms the solutions obtained by Gurobi using a 4-hour computational time budget (the average gap is -11.13% with the highest gap of 0.42% over the 20 test instances). The average CPU time for the metaheuristic is 31 seconds, showing its advantage for solving larger real-world applications.

5.4. Case study

Let us consider a typical application case where a feeder service in a disk service area connects to a central station. The station is located in the center and the depot is not far from the station. We also assume that a limited number of charging stations are available at the depot and at the station. Given a homogeneous fleet of EVs, we want to answer the following three questions.

    Q1: What is the impact of meeting point (MP) availability (parametrized by nearest separation distance) on system costs, service rate, and service level?
    Q2: What is the fleet size required to maintain a user-defined service rate given demand intensities?
    Q3: How does the algorithm perform for real-world problem sizes and demand intensities?

The first two experiments consider a circular service area with a 5 km radius. While the third one scales up the test instance to mimic a real-world example for the morning rush-hour commuting scenario from Arlon to Luxembourg. For the three experiments, the test instances are generated using our scenario generator described in Section 5.1. We assume that two charging stations are located at (0, 0) and (1000, 0), each with 2 DC fast chargers (50kWh). The feeder service provides 3 trips per hour with a regular headway of 20 minutes to connect to the transit station. Customers are randomly distributed in the ring between a radius of 1.5 km and the range limit of the service area. Figure 6 shows an example of the configuration of a test instance. A heterogeneous fleet is assumed with a vehicle capacity of 10 and 20 passengers, respectively. Energy consumption rate and battery capacity are shown in Table 1. To invoke potential charging operations for cases of real-world problem sizes, initial state of charge of vehicles are assumed to be randomly distributed between 30% and 80% at the beginning of the planning horizon. Vehicle speed is assumed as 30km/hr. All experiments are run on the same computer as described in the previous section.

    The demand scenarios and MP availability settings for the first two experiments are as follows. The demand distribution is assumed to follow a normal distribution concentrated between 7:40 and 8:20 (80%-90% demand) as shown on the right side of Figure 6. To analyze the impact of system parameters with different demand intensities, three demand levels are generated with 50, 100, and 200 customers (requests), respectively. Each dataset (corresponding to a demand level) contains 5 randomly generated test instances. A maximum walking distance of 1.0 km is assumed and the reference separation distance between two adjacent MPs is 1.2 km to generate a set of MPs. Note that an $X$ km MP separation distance results in the furthest away from the grid center of $X/\sqrt{2}$ km. So using 1.4 km as a MP separation distance results in $1.4/\sqrt{2} = 0.9899$ km farthest accessible MPs to customers.

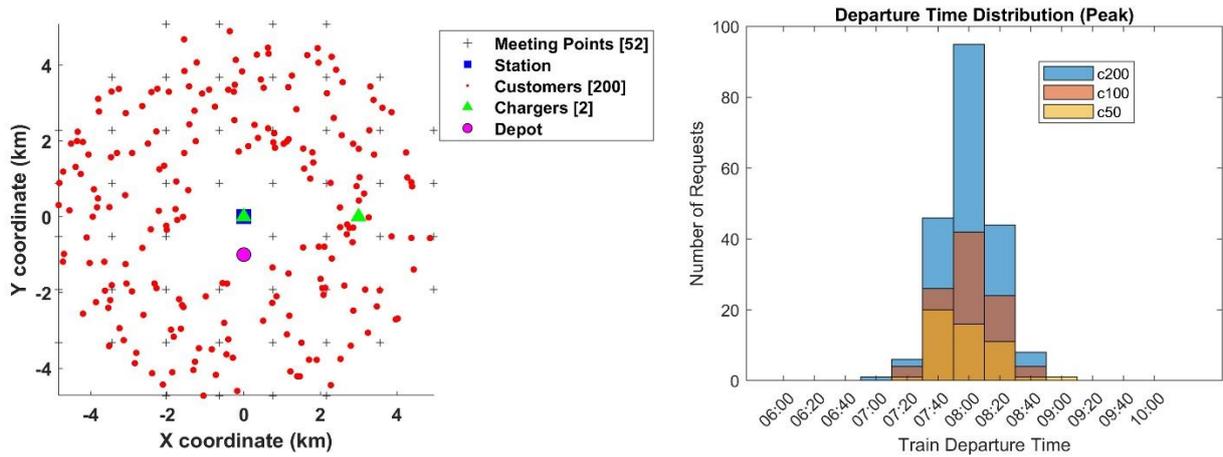

**Figure 6.** Example of a test instance with 200 customers on a disk area with a central train station at the center (left). Distributions of train departure times of requests (peak scenario) (right).

a. Experiment 1: Effect of meeting point separation distance

To analyze the impact of the distance between two nearest adjacent MPs, three levels of distance are considered, i.e., 1.0, 1.2, and 1.4 km. Using shorter separation distance can reduce the customer's walking distance at the cost of higher bus operating costs. On the other hand, using a distance greater than the maximum walking distance will reduce the service rate. A number of indicators are used to evaluate system costs and service levels of a solution: (a) number of vehicles used, (b) service rate: percentage of customers served, (c) walking distance: Euclidean distance between a customer's origins and their pick-up MPs, (d) in-vehicle travel time (IVT), (e) waiting time for buses arriving at the station ahead of schedule. (f) kilometers traveled by all vehicles in service (KMT), (g) number of customers served per kilometer traveled (cus/KMT), and (h) number of customers served per MP activated (cus/MP). All the reported indicators are collected based on the average where each test instance is solved 3 times by the metaheuristic to obtain the best solution.

Table 7 shows the system performance indicators for different demand levels and MP separation distances. Note that the fleet size for 50, 100 and 200 customers are (3, 2), (4, 5) and (6, 7) respectively, in which (n1, n2) means n1 10-seater buses and n2 20-seater buses. It can be observed that reducing the separation distance from 1.4 km to 1.0 km results in higher service rate (+~3% for c50 case (50 customers)), while it has little impact on the service rate for higher demands (c100 and c200 cases). Furthermore, the bus operator benefits from cost savings due to lower KMT (-4.9% on average for c100 and c200 cases, while +9.29% due to higher service rate). The value in brackets in the KMT column is the relative change with respect to the case with a separation distance of 1.4 km. The savings in KMT is due to better optimal customer-to-meeting-point assignment can be found when more MPs are available. If we look at the number of customer per vehicle KMT (Cus/KMT), higher demand contributes to a higher cost effectiveness (Cus/KMT is almost doubled when demand is increased from 50 to 200). Using a smaller MP separation distance- (1.0 km vs. 1.4 km) can increase Cus/KMT around 5% for both c100 and c200 cases. Note that in reality, demand locations are not randomly distributed in space, but rather concentrated in the city center or its suburbs. Therefore, these values can be amplified depending on the degree of demand concentration. The cus/MP column shows the average number of customers picked up together at the same stop. It can be observed that this number increases consistently when the MP separation distance reduces from 1.4 km to 1 km. The result quantifies the potential benefits for improving the system efficiency of the meeting-point-based on-demand feeder service.

Overall, we find that the operator can adjust this parameter to trade-off between operating costs and customer inconvenience by jointly considering the demand intensity and its spatial distribution to optimize its effect.

**Table 7. Sensitivity analysis to the availability of meeting points.**

| MP Dist* | Num cus | Num used veh | Service Rate (%) | Walking dist (km) | IVT (min) | KMT (km) | KMT (+-%) | cus/KMT | cus/MP | CPU (s) |
|---|---|---|---|---|---|---|---|---|---|---|
| 1.4 | 50  | 5    | 92.4  | 0.59 | 7.99 | 0       | 0       | 0.33 | 1.51 | 26 |
|     | 100 | 8.8  | 99.8  | 0.55 | 8.45 | 0       | 0       | 0.39 | 1.65 | 6  |
|     | 200 | 10.8 | 100.0 | 0.57 | 8.60 | 0       | 0       | 0.57 | 2.22 | 9  |
| 1.2 | 50  | 5    | 95.6  | 0.54 | 8.19 | 11.95%  | 11.95%  | 0.30 | 1.45 | 22 |
|     | 100 | 7.8  | 100.0 | 0.54 | 8.37 | -3.60%  | -3.60%  | 0.40 | 1.75 | 6  |
|     | 200 | 10.4 | 100.0 | 0.56 | 8.44 | 1.21%   | 1.21%   | 0.56 | 2.30 | 9  |
| 1.0 | 50  | 5    | 95.2  | 0.50 | 7.83 | 9.29%   | 9.29%   | 0.31 | 1.51 | 32 |
|     | 100 | 7.8  | 100.0 | 0.53 | 8.21 | -4.60%  | -4.60%  | 0.41 | 1.85 | 9  |
|     | 200 | 10.4 | 100.0 | 0.56 | 8.33 | -5.25%  | -5.25%  | 0.60 | 2.56 | 24 |

*MP dist: Distance between two closest adjacent MPs. The fleet sizes for the c50, c100, and c200 instances are 5, 9, and 14 vehicles, respectively.

b. Experiment 2: fleet size requirement

This experiment aims to answer a practical question: what is the trade-off between fleet size and service rate, given a user-defined meeting point separation distance? The previous experiment 1 shows that fleet size depends on both demand intensity and meeting point separation distance to meet customer demand. To answer this question, we consider the scenarios of 100 customers but vary the available fleet size from 3 to 9 over a vehicle pool of 5 10-seater buses and 4 20-seater buses). The MP separation distance is also varied from 1.4 km to 1.0 km to analyze their effect on the service rate and the operating cost in terms of the number of kilometers traveled by the fleet. $E_{init}^k$ is set as 80% of their respective battery capacities to isolate the impact of these two factors. Note that our algorithm sorts the driving range (i.e. initial state of charge of vehicles divided by their unitary energy consumption) of vehicles in descending order to select vehicles to use. In this way, vehicle charging needs can be always reduced.

The left part of Figure 7 shows that the service rate increases faster when the fleet size is no more than 6 vehicles. The service rate reaches 95.2% and 96.0% when shorter separation distances (1.2 km and 1.0 km, respectively) are used. When a larger separation distance (1.4 km) is used, the service rate drops to 91% (-5.0%). Fleet mileage shows similar curves with increasing mileage when a larger fleet is used (right side of Figure 7). However, the fleet mileage is reduced when lower MP separation distances are used. The results confirm that using lower MP separation distance is more beneficial to optimize the system performance (higher service rate and more cost effective). Thus, the range for setting the fleet size is between 6 and 8 vehicles to satisfy at least 95% of the demand with at least the same maximum walking distance of customers.

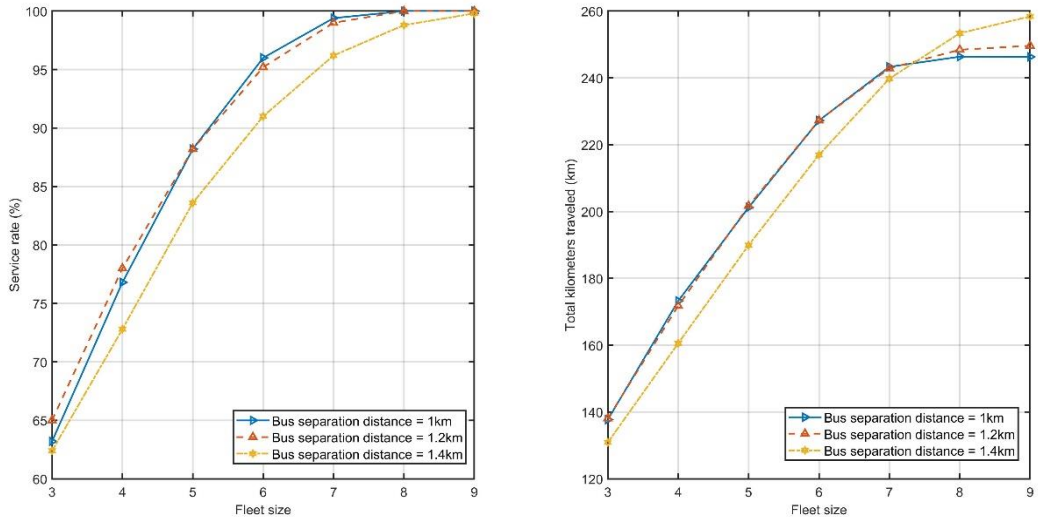

**Figure 7. The effect of fleet size on customer service rates and fleet mileage for different levels of meeting point separation distances (number of customers = 100).**

c. Experiment 3: Larger test instances

We consider a real larger case in the Belgian side of the border area of Luxembourg. According to recent statistics, Arlon accounts for 9% of Luxembourg's cross border workers in 2023, with around 16900 workers. The analysis of the 2017 Luxmobil survey data shows that the Belgium-side cross border workers mainly use the car for their daily commute (88.2%). This is followed by public transport (11.8%), where 7.9% of all trips are made by train (Lambotte et al., 2021, p. 27). Thus, a total of 1335 train commuting trips per day can be estimated from and around Arlon to Luxembourg. The cross-border region of Arlon and the rail network connection to Luxembourg is shown in Figure 8. In summary, there is only one rail line (line 50) between Arlon and Luxembourg. During the rush hour (6:00-9:00), there are 13 trains from Arlon to Luxembourg, i.e. an average frequency of 15 minutes, according to the current timetables of this line. Arlon covers an area of 118.6 km2 with the central station in the city center (Figure 8). This corresponds to a circle with a radius of 6.1442 km. Based on this real-world application case, we use our scenario generator to generate a test instance with a radius of 6 km and assume that customers are located within a ring of 1.5 km and 6 km from the center. We assume that the maximum walking distance of customers is 1 km. A meeting point separation distance of 0.8 km is used to generate the locations of MPs. This results in a total of 111 MPs. We assume that the operator provides three services per hour (i.e. buses arrive at the station every 20 minutes, with a buffer time of 10 minutes) The number of requests is 1000 (75% of the 1335 train commuting trips), randomly distributed in the service area and in the desired arrival times at the station (i.e. 6:00, 6:20, 6:40, 7:00, 7:20, and 7:40). The characteristics of the fleet, the average speed of the vehicles, the locations of the depot, the station, and the charging stations are the same as in the two previous experiments. Initial state of charge of vehicles are set randomly between 30%-80% of their battery capacities (58.28% on average).

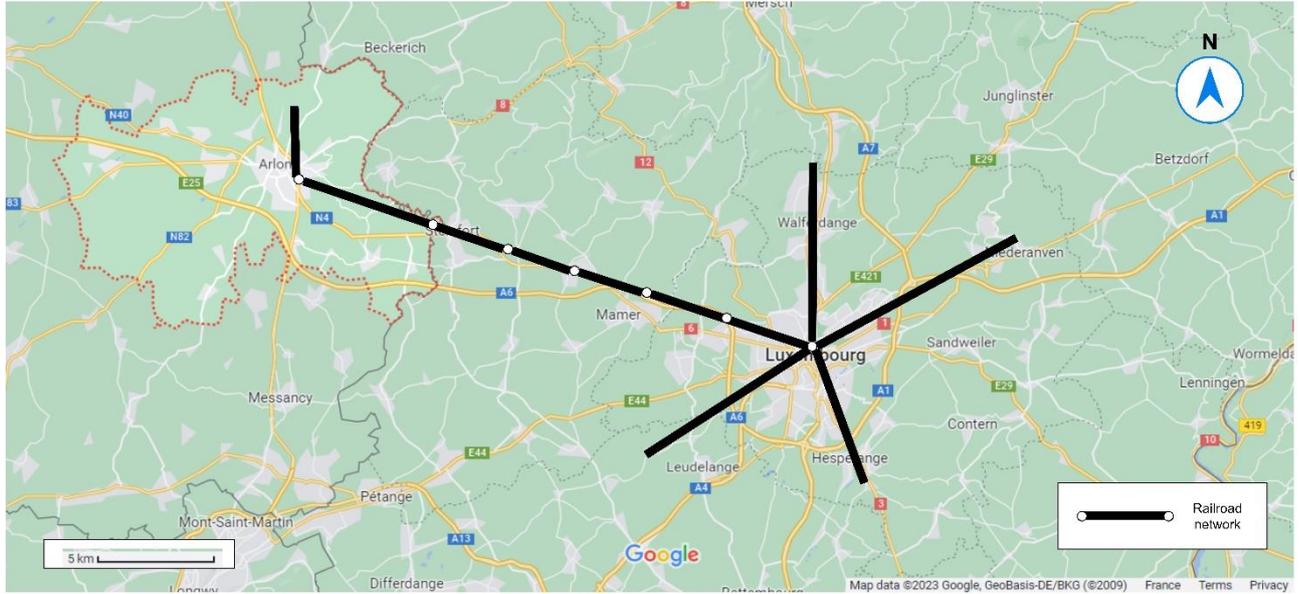

**Figure 8. The region of Arlon and the railway connection to Luxembourg.**

In summary, the problem size of this instance after preprocessing is |R|=1000, |V|=653, |K|=18, $|\mathcal{A}_B|$=77828, and $|\mathcal{A}_C|$=2830, respectively. The fleet size varies from 10 to 18 with 50% of each vehicle type to evaluate its impact on the system performance. As shown in Table 8, we can solve this instance using the hybrid metaheuristic in about 15 minutes for the hardest case with 14 vehicles. The charging time ranges from 10.7 minutes to 21 minutes. The service rate is about 68.7% when using 10 vehicles, while it is 98.6% when using 18 vehicles. The average walking distance is about 0.6 km, while the average vehicle travel time is about 10 minutes. The average number of customers assigned to an activated stop is about 4, and the total KMT depends on the number of served customers. Notice that for solving each MILP problem in the post-optimization procedure (Algorithm 3), a 120-second time limit is set to load the model and obtain a feasible solution.

**Table 8. Performance of the feeder service with respect to different fleet sizes (number of customers=1000).**

| |K| | | Service rate (%) | cus unserved | Walking dist (km) | IVT (min) | KMT (km) | CT* (min) | cus/KMT | cus/MP | CPU (s) |
|---|---|---|---|---|---|---|---|---|---|---|
| T1 bus* | T2 bus* | | | | | | | | | |
| 9 | 9 | 98.6 | 14 | 0.57 | 9.9 | 1030.2 | 10.7 | 0.96 | 4.11 | 473 |
| 8 | 8 | 97.6 | 24 | 0.58 | 9.9 | 997.9 | 21.0 | 0.98 | 4.21 | 608 |
| 7 | 7 | 85.7 | 143 | 0.58 | 9.6 | 851.5 | 10.9 | 1.01 | 4.10 | 877 |
| 6 | 6 | 79.2 | 208 | 0.58 | 9.6 | 739.7 | 12.0 | 1.07 | 3.98 | 676 |
| 5 | 5 | 68.7 | 313 | 0.57 | 9.3 | 624.4 | 16.4 | 1.10 | 4.07 | 636 |

Remark: 1. T1 (T2): Type-1 (Type-2) bus, corresponding to the two bus types (10- and 20- seaters) in Table 1; 2. CT: Total charging time.

When further analyzing the impact of the post-optimization procedure, we found that it allows to serve more customers from +0.1% to +4.2% (42/1000 customers) compared with that without post-optimization. The problem size of the MILP in the matheuristic and the number of times that it has to be solved ($|\mathcal{L}_{\bar{R}}|$) are shown in Table 9.

Table 9. Comparison of the performance of the algorithm with and without post-optimization.

| Num. of vehicles | Without post-optimization | | | With post-optimization | | | | | | | |
|---|---|---|---|---|---|---|---|---|---|---|---|
| | Obj. value | Service rate (%) | CPU (s) | Obj. value | Service rate (%) | CPU (s) | Problem size of the MILP in the post-optimization* | | | | |
| | | | | | | | $|R_\ell|$ | $|K_\ell|$ | $|G_\ell|$ | $|\mathcal{A}_\ell|$ | $|\mathcal{L}_{\tilde{R}}|$ |
| 18 | 9375.0 | 98.5 | 381 | 9202.2 | 98.6 (+0.1) | 473 | 167 | 15 | 107 | 14552 | 1 |
| 16 | 10302.0 | 95.5 | 373 | 9671.8 | 97.6 (+2.1) | 608 | 170 | 15 | 107 | 14459 | 5 |
| 14 | 14329.6 | 82.6 | 180 | 13291.9 | 85.7 (+3.1) | 877 | 167 | 13 | 106 | 13968 | 6 |
| 12 | 16676.0 | 75.0 | 156 | 15180.9 | 79.2 (+4.2) | 676 | 167 | 12 | 108 | 13960 | 6 |
| 10 | 19855.2 | 64.5 | 198 | 18428.9 | 68.7 (+4.2) | 636 | 167 | 10 | 107 | 13407 | 6 |

*$|X|$: number of elements in *X*.

## 6. Conclusions and discussion

Electric vehicle routing and charging scheduling problems are of great importance in passenger transportation and logistics applications, and have received increasing research interest for decades. For passenger transportation, the problem is closely related to the dial-a-ride problem and its variants using electric vehicles. Previous studies have shown that the classical door-to-door based DRT system can be more cost-effective if customers can be picked up or dropped off at predefined meeting points. Several recent studies and real-world applications have demonstrated the benefits of an on-demand meeting point-based DRT system to achieve significant cost savings. However, when deploying such a system using a fleet of electric vehicles, the resulting routing problem needs to jointly consider vehicle charging scheduling with partial charging capability given limited charging station capacity constraints. However, the above routing problems with charging synchronization are more difficult to solve, and efficient solution algorithms for solving medium/large problem instances are still underdeveloped.

In this study, we considered the problem of an electric on-demand meeting-point-based feeder system with charging synchronization constraints and proposed a layered graph model and a mixture of randomization and greedy strategy within a hybrid DA-based algorithm framework to solve this problem efficiently. The originality of this paper is to jointly optimize customer-to-meeting-point assignment, electric bus routing and charging scheduling. The new layered graph structure is proposed to trim unnecessary nodes and arcs to reduce the problem size. We customized the DA algorithm using two randomization methods, namely relocate ensemble and destroy and repair, to improve the exploration of the solution space. Moreover, for the customer-to-meeting-point assignment, instead of using the same parameter (rho) to trade off users' walking time and bus routing time, we allow a layer-specific rho to be tuned when there are unserved customers for certain layers. In this way, the tuned parameter is more reactive to demand intensity variations.

We tested the algorithm on 20 test instances with up to 100 customers and 49 meeting points under different initial battery levels and demand distributions (peak and off-peak scenarios) and compare it with an exact solution method. The results show that the proposed algorithm can efficiently find solutions with good solution quality. The performance of the proposed algorithm is compared with the solutions obtained by a MILP solver given a computation time of 4 hours. The results show that the metaheuristic can obtain good solutions (average and best gap of 1.80% and 1.69%, respectively) within a short computation time (less than 1 minute on average). Sensitivity analysis is performed to investigate the effect of the algorithmic parameters. To stimulate real-world applications, a case study is conducted to analyze the effect of the system parameters (meeting point separation distance and fleet size) and to solve a real-world problem size case with 1000 customers and 111 meeting points.

Several research directions are ongoing, including joint planning of charging infrastructure and fleet size, and integrated optimization of DRT system operating policies. Future extensions could consider its applications to real-world case studies, extending the current static model to a dynamic one to allow the insertion of new customers. Regarding the modeling of charging behavior, more realistic charging functions and can be integrated in the future (non-linear charging functions (Froger et al., 2022)). In addition, it would be beneficial to consider soft (allowing vehicles to wait at charging stations) instead of hard charging station capacity constraints, and to introduce more realistic components when modeling the vehicle recharging process, such as introducing a fixed cost/time for each recharge, among others, and time-dependent waiting times when public recharging stations are involved (Keskin et al., 2019). Furthermore, there is still room for improvements in the performance of the proposed metaheuristic (e.g. employing and testing different methods such as genetic algorithm), which is our ongoing research.

**Acknowledgements:** The work was supported by the Luxembourg National Research Fund (C20/SC/14703944).

**References**

1. Aïvodji, U.M., Gambs, S., Huguet, M.-J., Killijian, M.-O., 2016. Meeting points in ridesharing: A privacy-preserving approach. Transp. Res. Part C Emerg. Technol. 72, 239–253. https://doi.org/10.1016/j.trc.2016.09.017
2. Alonso-González, M.J., Liu, T., Cats, O., Van Oort, N., Hoogendoorn, S., 2018. The Potential of Demand-Responsive Transport as a Complement to Public Transport: An Assessment Framework and an Empirical Evaluation. Transp. Res. Rec. 2672. https://doi.org/10.1177/0361198118790842
3. Arnold, F., Sörensen, K., 2019. Knowledge-guided local search for the vehicle routing problem. Comput. Oper. Res. 105, 32–46. https://doi.org/10.1016/j.cor.2019.01.002
4. Belenguer, J.M., Benavent, E., Prins, C., Prodhon, C., Wolfler Calvo, R., 2011. A Branch-and-Cut method for the Capacitated Location-Routing Problem. Comput. Oper. Res. 38, 931–941. https://doi.org/10.1016/J.COR.2010.09.019
5. Bian, Z., Liu, X., 2019. Mechanism design for first-mile ridesharing based on personalized requirements part I: Theoretical analysis in generalized scenarios. Transp. Res. Part B Methodol. 120, 147–171. https://doi.org/10.1016/j.trb.2018.12.009
6. Bongiovanni, C., Kaspi, M., Geroliminis, N., 2019. The electric autonomous dial-a-ride problem. Transp. Res. Part B Methodol. 122, 436–456. https://doi.org/10.1016/j.trb.2019.03.004
7. Braekers, K., Caris, A., Janssens, G.K., 2014. Exact and meta-heuristic approach for a general heterogeneous dial-a-ride problem with multiple depots. Transp. Res. Part B Methodol. 67, 166–186. https://doi.org/10.1016/j.trb.2014.05.007
8. Bruglieri, M., Mancini, S., Pisacane, O., 2021. A more efficient cutting planes approach for the green vehicle routing problem with capacitated alternative fuel stations. Optim. Lett. 15, 2813–2829. https://doi.org/10.1007/S11590-021-01714-3/FIGURES/3
9. Bruglieri, M., Mancini, S., Pisacane, O., 2019. The green vehicle routing problem with capacitated alternative fuel stations. Comput. Oper. Res. 112. https://doi.org/10.1016/j.cor.2019.07.017
10. Chen, P.W., Nie, Y.M., 2017. Analysis of an idealized system of demand adaptive paired-line hybrid transit. Transp. Res. Part B Methodol. 102, 38–54. https://doi.org/10.1016/j.trb.2017.05.004
11. Cordeau, J.F., Laporte, G., 2003. A tabu search heuristic for the static multi-vehicle dial-a-ride problem. Transp. Res. Part B Methodol. 37, 579–594. https://doi.org/10.1016/S0191-2615(02)00045-0
12. Czioska, P., Kutadinata, R., Trifunović, A., Winter, S., Sester, M., Friedrich, B., 2019. Real-world meeting points for shared demand-responsive transportation systems, Public Transport. Springer Berlin Heidelberg. https://doi.org/10.1007/s12469-019-00207-y
13. Desaulniers, G., Errico, F., Irnich, S., Schneider, M., 2016. Exact algorithms for electric vehicle-routing problems with time windows. Oper. Res. 64, 1388–1405. https://doi.org/10.1287/opre.2016.1535
14. Dueck, G., Scheuer, T. 1990. Threshold accepting: A general purpose optimization algorithm appearing superior to simulated annealing. Journal of computational physics, 90(1), 161-175.
15. Erdoğan, S., Miller-Hooks, E., 2012. A Green Vehicle Routing Problem. Transp. Res. Part E Logist. Transp. Rev. 48, 100–114. https://doi.org/10.1016/j.tre.2011.08.001
16. Felipe, Á., Ortuño, M.T., Righini, G., Tirado, G., 2014. A heuristic approach for the green vehicle routing problem with multiple technologies and partial recharges. Transp. Res. Part E Logist. Transp. Rev. 71, 111–128. https://doi.org/10.1016/j.tre.2014.09.003

17. Froger, A., Jabali, O., Mendoza, J.E., Laporte, G., 2021. The Electric Vehicle Routing Problem with Capacitated Charging Stations. Transp. Sci. 56, 460–482. https://doi.org/10.1287/TRSC.2021.1111
18. Froger, A., Mendoza, J.E., Jabali, O., Laporte, G., 2017. A Matheuristic for the Electric Vehicle Routing Problem with Capacitated Charging Stations. Working paper. CIRRELT-2017-31.
19. Goeke, D., Schneider, M., 2015. Routing a mixed fleet of electric and conventional vehicles. Eur. J. Oper. Res. 245, 81–99. https://doi.org/10.1016/j.ejor.2015.01.049
20. Haglund, N., Mladenović, M.N., Kujala, R., Weckström, C., Saramäki, J., 2019. Where did Kutsuplus drive us? Ex post evaluation of on-demand micro-transit pilot in the Helsinki capital region. Res. Transp. Bus. Manag. 32. https://doi.org/10.1016/j.rtbm.2019.100390
21. Jenn, A. (2019). Electrifying Ride-Sharing: Transitioning to a Cleaner Future. UC Davis: National Center for Sustainable Transportation. Retrieved from https://escholarship.org/uc/item/12s554kd
22. Keskin, M., Çatay, B., 2016. Partial recharge strategies for the electric vehicle routing problem with time windows. Transp. Res. Part C Emerg. Technol. 65, 111–127. https://doi.org/10.1016/j.trc.2016.01.013
23. Keskin, M., Laporte, G., Çatay, B., 2019. Electric Vehicle Routing Problem with Time-Dependent Waiting Times at Recharging Stations. Comput. Oper. Res. 107, 77–94. https://doi.org/10.1016/j.cor.2019.02.014
24. Kucukoglu, I., Dewil, R., Cattrysse, D., 2021. The electric vehicle routing problem and its variations: A literature review. Comput. Ind. Eng. 161, 107650. https://doi.org/10.1016/j.cie.2021.107650
25. Lam, E., Desaulniers, G., Stuckey, P.J., 2022. Branch-and-cut-and-price for the Electric Vehicle Routing Problem with Time Windows, Piecewise-Linear Recharging and Capacitated Recharging Stations. Comput. Oper. Res. 145, 105870. https://doi.org/10.1016/J.COR.2022.105870
26. Lambotte, J. M., Marbehant, S., Rouchet, H., 2021. Distribution spatiale de l'emploi et des modes de transport des travailleurs actifs au Grand-Duché de Luxembourg–Analyse des données de l'enquête Luxmobil 2017. UniGR-CBS Working Paper, 11.
27. Lee, A., Savelsbergh, M., 2017. An extended demand responsive connector. EURO J. Transp. Logist. 6, 25–50. https://doi.org/10.1007/s13676-014-0060-6
28. Lutz, R. 2004. Adaptive Large Neighborhood Search. Bachelor thesis at Ulm University.
29. Ma, T.-Y., Rasulkhani, S., Chow, J.Y.J., Klein, S., 2019. A dynamic ridesharing dispatch and idle vehicle repositioning strategy with integrated transit transfers. Transp. Res. Part E Logist. Transp. Rev. 128, 417–442. https://doi.org/10.1016/j.tre.2019.07.002
30. Ma, T.Y., Chow, J.Y.J., Klein, S., Ma, Z., 2021. A user-operator assignment game with heterogeneous user groups for empirical evaluation of a microtransit service in Luxembourg. Transp. A Transp. Sci. 17, 946–973. https://doi.org/10.1080/23249935.2020.1820625
31. Macrina, G., Di Puglia Pugliese, L., Guerriero, F., Laporte, G., 2019. The green mixed fleet vehicle routing problem with partial battery recharging and time windows. Comput. Oper. Res. 101, 183–199. https://doi.org/10.1016/j.cor.2018.07.012
32. Melis, L., Queiroz, M., Sörensen, K., 2021. The integrated on-demand bus routing problem. Working Papers 2021004, University of Antwerp, Faculty of Business and Economics.
33. Melis, L., Sörensen, K., 2022. The static on-demand bus routing problem: large neighborhood search for a dial-a-ride problem with bus station assignment. Int. Trans. Oper. Res. 29, 1417–1453. https://doi.org/10.1111/itor.13058
34. Montenegro, B.D.G., Sörensen, K., Vansteenwegen, P., 2021. A large neighborhood search algorithm to optimize a demand-responsive feeder service. Transp. Res. Part C Emerg. Technol. 127, 103102. https://doi.org/10.1016/j.trc.2021.103102
35. Montenegro, B.D.G., Sörensen, K., Vansteenwegen, P., 2022. A column generation algorithm for the demand-responsive feeder service with mandatory and optional, clustered bus-stops. Networks 80, 274–296. https://doi.org/10.1002/net.22095
36. Park, J., Kim, B.-I., 2010. The school bus routing problem: A review. Eur. J. Oper. Res. 202, 311–319. https://doi.org/10.1016/j.ejor.2009.05.017
37. Parragh, S.N., Doerner, K.F., Hartl, R.F., 2010. Variable neighborhood search for the dial-a-ride problem. Comput. Oper. Res. 37, 1129–1138. https://doi.org/10.1016/j.cor.2009.10.003

38. Pimenta, V., Quilliot, A., Toussaint, H., Vigo, D., 2017. Models and algorithms for reliability-oriented Dial-a-Ride with autonomous electric vehicles. Eur. J. Oper. Res. 257, 601–613. https://doi.org/10.1016/j.ejor.2016.07.037
39. Ropke, S., Pisinger, D., 2006. An adaptive large neighborhood search heuristic for the pickup and delivery problem with time windows. Transp. Sci. 40, 455–472. https://doi.org/10.1287/trsc.1050.0135
40. Schittekat, P., Kinable, J., Sörensen, K., Sevaux, M., Spieksma, F., Springael, J., 2013. A metaheuristic for the school bus routing problem with bus stop selection. Eur. J. Oper. Res. 229, 518–528. https://doi.org/10.1016/j.ejor.2013.02.025
41. Schneider, M., Stenger, A., Goeke, D., 2014. The electric vehicle-routing problem with time windows and recharging stations. Transp. Sci. 48, 500–520. https://doi.org/10.1287/trsc.2013.0490
42. Stiglic, M., Agatz, N., Savelsbergh, M., Gradisar, M., 2015. The benefits of meeting points in ride-sharing systems. Transp. Res. Part B Methodol. 82, 36–53. https://doi.org/10.1016/j.trb.2015.07.025
43. Vansteenwegen, P., Melis, L., Aktaş, D., Montenegro, B.D.G., Sartori Vieira, F., Sörensen, K., 2022. A survey on demand-responsive public bus systems. Transp. Res. Part C Emerg. Technol. 137, 103573. https://doi.org/10.1016/j.trc.2022.103573
44. Wang, F., Ross, C.L., 2019. New potential for multimodal connection: exploring the relationship between taxi and transit in New York City (NYC). Transportation (Amst). 46, 1051–1072. https://doi.org/10.1007/s11116-017-9787-x
45. Zheng, Y., Li, W., Qiu, F., Wei, H., 2019. The benefits of introducing meeting points into flex-route transit services. Transp. Res. Part C Emerg. Technol. 106, 98–112. https://doi.org/10.1016/j.trc.2019.07.012

**Appendix A. Notation**

| | |
|---|---|
| *Sets* | |
| $G$ | Set of physical meeting points, i.e. $G = \{1, \ldots, N_G\}$ |
| $G'$ | Set of dummy (duplicate) meeting point vertices (nodes) |
| $D$ | Set of physical transit stations, i.e. $D = \{1, \ldots, N_D\}$ |
| $D'$ | Set of dummy (duplicate) transit station vertices |
| $S$ | Set of physical chargers, i.e. $S = \{1, \ldots, N_S\}$ |
| $S'$ | Set of dummy (duplicate) charger vertices |
| $R$ | Set of customers (i.e. location of origin of customers) |
| $K$ | Set of electric buses |
| $\bar{V}$ | Set of all vertices, i.e. $\bar{V} = G' \cup D' \cup S' \cup R \cup \{0, N+1\}$ |
| $V$ | Subset of vertices, i.e. $V = G' \cup D' \cup S'$ |
| $V_0, V_{N+1}, V_{0,N+1}$ | $V_0 = V \cup \{0\}, V_{N+1} = V \cup \{N+1\}, V_{0,N+1} = V \cup \{0, N+1\}$ |
| $\mathcal{A}_C$ | Set of walking arcs from customers' origins to meeting points, i.e. $\mathcal{A}_C = \{(r, j) \mid r \in R, j \in G'\}$ |
| $\mathcal{A}_B$ | Set of bus arcs |
| *Parameters and auxiliary variables* | |
| $0, N+1$ | Two duplicate instances of the depot |
| $T$ | Planning horizon |
| $A_i^k$ | Arrival time of bus $k$ at vertex $i$ |
| $B_i^k$ | Beginning time of service of bus $k$ at vertex $i$ |
| $p_i^k$ | Indicator: 1 if node $i$ is visited by bus k, 0 otherwise |
| $W_i^k$ | Waiting time of bus $k$ at node $i$ |
| $q_i^k$ | Passenger load of bus $k$ when leaving vertex $i$ |
| $E_i^k$ | Battery energy level of bus $k$ when arriving at the vertex $i$ |
| $Q^k$ | Capacity of bus $k$ |
| $E_{min}^k, E_{max}^k, E_{init}^k$ | Minimum, maximum, initial state of charge (SOC) of bus $k$ |
| $w_{ri}$ | Walking distance from customer $r$ origin to meeting point $i \in G'$ |
| $c_{ij}$ | Distance from vertex $i$ to vertex $j$ |
| $t_{ij}$ | Bus travel time from vertex $i$ to vertex $j$ $\forall i, j \in V_{0,N+1}$. Note that $t_{rj}$ is the walking time from customer $r$ origin to meeting point $j$, $\forall r \in R, j \in G'$. |
| $L_i$ | Maximum ride time for customers picked-up at node $i$. Calculated as 'straight line' ride time multiplied by a detour factor. |
| $w_{max}$ | Maximum walking distance |
| $u_i$ | Service time at vertex $i \in V$ |
| $e_i, l_i$ | Earliest and latest starting times of service at vertex $i \in V_{0,N+1}$ |
| $d_r, d_i$ | Drop-off transit station dummy node of customer $r \in R$ and $i \in G'$ |
| $\omega$ | Penalty cost if a customer is rejected |
| $\alpha_s$ | Charging rate of charger $s \in S'$ |
| $\beta^k$ | Energy consumption rate per kilometer traveled for bus $k$ |
| $M$ | Large positive number |
| $\lambda$ | Weighting coefficient for the objective function |
| $\rho$ | Parameter used in the customer-to-meeting-point assignment problem (Eq. (36)) |
| $\theta_i$ | Indicator being 1 if a dummy node is visited by a bus (used in Eqs. (36)-(46)) |
| *Decision variables* | |
| $y_{ri}^k$ | Indicator: 1 if customer $r$ is assigned to bus $k$ and meeting point $i$, 0 otherwise |
| $x_{ij}^k$ | Indicator: 1 if arc $(i, j)$ is traversed by bus $k$, 0 otherwise |
| $\tau_s^k$ | Charging duration for bus $k$ at charger $s$, $s \in S'$ |

## Appendix B. Algorithm for reassigning unserved customers

This algorithm aims at reassigning the customers (including both served and unserved customers) of certain layers containing unserved customers to other meeting points in order to improve the current solution and increase the bus service rate. Let $s$ and $s'$ denote the current and temporary solution (both are initialized as the solution obtained by DA algorithm), respectively. The algorithm first obtains the set of unserved customers of $s$, denoted as $\tilde{R}$, and their corresponding layers $\mathcal{L}_{\tilde{R}}$. Let $\tilde{R}_\ell$ and $R_\ell$ be the subsets of the unserved customers and all customers on layer $\ell$, respectively. For each of these layers, the algorithm solves a customer reassignment problem formulated as a MILP of Eqs. (B1)-(B16) with a short computation time limit (e.g., tens of seconds depending on the problem size). This MILP aims to find an alternative meeting point assignment and partial routes under the constraints of vehicle capacity, maximum walking distance, and maximum ride time of customers (line 5 in Algorithm 3). Note that the subroutes on the layers without unserved customers remain unchanged. After getting the solutions of the above MILP problems, the time window constraints are checked (line 6 in Algorithm 3). If these constraints are satisfied, the preliminary solution is updated. If some routes of $s'$ have been modified, check the constraints with respect to the state of charge of these routes (line 10 of Algorithm 3). If violated, schedule vehicle charging operations in a sequential manner using a greedy charging rescheduling policy without conflicts with the existing charging schedules of other vehicles (lines 12-20 of Algorithm 3). Finally, if the temporary solution has a lower objective function value, return $s'$. Figure B1 shows an illustrative example of customer reassignment on a layer with unserved customers. Given the penalty associated with unserved customers, the obtained solution tends to accommodate as many customers as possible.

Algorithm 3. Customer reassignment algorithm for the layers with unserved customers.

1. **Input**: solution $s$
2. Let $\tilde{R} :=$ set of unserved customers of $s$, $\mathcal{L}_{\tilde{R}} :=$ set of layers contain $\tilde{R}$, $R_{\mathcal{L}_{\tilde{R}}} :=$ set of all customers on $L_{\tilde{R}}$, $K_\ell :=$ set of vehicles traversing layer $\ell$.
3. **for** $\ell \in \mathcal{L}_{\tilde{R}}$
4.    identify the precedent nodes and the transit nodes for $k \in K_\ell$ traversing layer $\ell$.
5.    solve the MILP of Eqs. (B1)-(B16) within the user-defined computation time limit.
6.    **If** the optimal solution is found and the time window constraints for the obtained partial routes are satisfied,
7.      save the temporary solutions (partial routes of vehicles traversing layer $\ell$).
8.    **end**
9. **end**
10. Let $K_1$ denote the set of vehicles with newly changed routes, and $K_2$ denotes the set of vehicles with unchanged routes. Create a charging state occupancy lookup table for all chargers based on the charging schedules of the vehicles in $K_2$. Let $success = true$
11. **for** $k \in K_1$
12.    **if** the state of charge of $k$ is insufficient, apply the greedy charging rescheduling algorithm without violating
13.    the capacity of each charger.
14.      **if** the charging scheduling is successful
15.        update the occupancy lookup table with newly added charging events of $k$.
16.      **else**
17.        $success = false$; return
18.      **end**
19.    **end**
20. **end**
21. Let $s^*$ denote the new temporary solution.
22. **if** $success = true$ and $c(s^*) < c(s)$
23.    return $s^*$
24. **end**

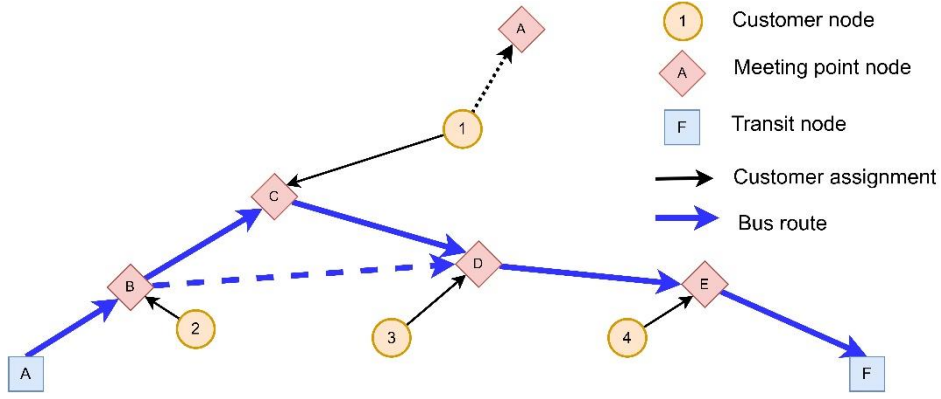

**Figure B1. A partial bus route from A to F after inserting unserved customer 1. The solid lines are the bus route or customer-to-meeting-point reassignments and the dashed lines are the removed bus routes and removed customer assignments.**

Let $K_\ell$ denote the set of vehicles traversing layer $\ell$. Consider a layer $\ell \in \mathcal{L}_{\tilde{R}}$, the subroutes to be changed are the visit sequences of meeting point nodes of $K_\ell$ on layer $\ell$. Let $o_\ell^k$ denote the predecessor node of the first meeting point on the subroute of vehicle k on layer $\ell$. Let $d_\ell^k$ denote the dummy transit station node on layer $\ell$ visited by vehicle $k$. The problem is a variant of the constrained vehicle routing problem with multiple vehicles starting at $o_\ell^k$ and ending at $d_\ell^k$ for $\forall k \in K_\ell$. The MILP formulation for reassigning customers is as follows.

$$\text{Min} \sum_{k \in K_\ell} \sum_{(i,j) \in \mathcal{A}_\ell} t_{ij} x_{ij}^k + \omega \sum_{r \in R_\ell} \left(1 - \sum_{i \in G'_\ell} y_{ri}\right) \tag{B1}$$

$$\sum_{i \in G'_\ell} y_{ri} \le 1, \forall r \in R_\ell \tag{B2}$$

$$t_{ri} y_{ri} \le w_{max}, \forall r \in R_\ell, i \in G'_\ell \tag{B3}$$

$$\sum_{r \in R_\ell} y_{ri} \le M\theta_i, \forall i \in G'_\ell \tag{B4}$$

$$\sum_{i \in G'_\ell} x_{o_\ell^k i}^k = 1, \forall k \in K_\ell \tag{B5}$$

$$\sum_{i \in G'_\ell} x_{id_\ell^k}^k = 1, \forall k \in K_\ell \tag{B6}$$

$$\sum_{i \in G'_\ell \cup o_\ell^k} x_{ij}^k - \sum_{i \in G'_\ell \cup d_\ell^k} x_{ji}^k = 0, \forall k \in K_\ell, j \in G'_\ell \tag{B7}$$

$$\theta_i = \sum_{k \in K_\ell} \sum_{j \in G'_\ell \cup d_\ell^k} x_{ij}^k, \forall i \in G'_\ell \tag{B8}$$

$$A_{o_\ell^k}^k = 0, \forall k \in K_\ell \tag{B9}$$

$$A_j^k \ge A_i^k + u_i + t_{ij} - M_2(1 - x_{ij}^k), \forall k \in K_\ell, (i,j) \in A_\ell \tag{B10}$$

$$A_{d_\ell^k}^k - A_i^k - u_i \le L_i + M_2(1 - x_{ij}^k), \forall k \in K_\ell, (i,j) \in A_\ell \tag{B11}$$

$$x_{ij}^k = 1 \Rightarrow q_j^k = q_i^k + \sum_{r \in R_\ell} y_{rj}, \forall k \in K_\ell, (i,j) \in A_{G'_\ell} \tag{B12}$$

$$0 \le q_i^k \le Q^k, \quad \forall k \in K_\ell, i \in G'_\ell \tag{B13}$$

$$q^k_{o^k_\ell} = 0, \forall k \in K_\ell \tag{B14}$$

$$y_{ri} \in \{0,1\}, \theta_i \in \{0,1\}, \forall r \in R_\ell, i \in G'_\ell \tag{B15}$$

$$A^k_i \geq 0, i \in G'_\ell \cup o^k_\ell \cup d^k_\ell, k \in K_\ell \tag{B16}$$

Given a layer $\ell \in \mathcal{L}_{\tilde{R}}$, the objective function (B1) minimizes the total bus routing times and the total penalties of unserved customers on layer $\ell$. Note that when attempting to serve more customers in the post-optimization procedure, customers' walking time minimization is relaxed to allow better re-organizing customer-to-meeting-point assignment in favor of serving more customers. The cost of such relaxation is recompensed by reducing the total penalty of unserved customers as the penalty of an unserved customer is relatively high. Constraints (B2)-(B3) state that each customer can be assigned to at most one meeting point within maximum walking distance. Constraints (B4) and (B8) state that each assigned meeting point must be visited by exactly one vehicle. Constraints (B5)-(B7) are the vehicle flow conservation constraints. Constraint (B9) sets the starting time of vehicles to 0 at $o^k_\ell$ for $\forall k \in K_\ell$. Constraint (B10) specifies the arrival times of vehicles at the sequence of visited meeting-point/transit nodes, while constraint (B11) specifies the riding time constraint of customers. Constraints (B12)-(B14) ensure that the vehicle load cannot exceed its capacity. Note that the initial load is 0 because the previous customers have been dropped off. $M = |R_\ell|$ and $\mathcal{A}_\ell$ denote the set of bus routing arcs on layer $\ell$. The notation of the other variables is the same as described in Appendix A.

**Appendix C. Comparison of the performance of the hybrid metaheuristic with and without the layered graph structure**

To show the effectiveness of the layered graph structure, we test the performance of the hybrid metaheuristic algorithm with and without using this structure in the customer-to-meeting-point assignment. To be clear, the MILP for the customer-to-meeting-point assignment problem without using the layered graph structure is shown below. We can observe that without using the layered graph structure, the solution spaces are much larger.

$$\text{Min } \lambda_2 \sum_{r \in R} \sum_{j \in G'} t_{rj} y_{rj} + \lambda_1 \rho \sum_{i \in G'} \sum_{j \in G' | \ell(i) = \ell(j)} t_{ij} z_{ij} \tag{C1}$$

Subject to
$$c_{rj} y_{rj} \leq w_{max}, \forall r \in R, j \in G' \tag{C2}$$

$$\sum_{j \in G'} y_{rj} = 1, \forall r \in R \tag{C3}$$

$$\sum_{r \in R} y_{rj} \leq Q_{max}, \forall j \in G' \tag{C4}$$

$$\sum_{r \in R} y_{rj} \leq M \theta_j, \forall j \in G' \tag{C5}$$

$$z_{ij} \leq \theta_i, \forall i \in G', j \in G' \tag{C6}$$

$$z_{ij} \leq \theta_j, \forall i \in G', j \in G' \tag{C7}$$

$$z_{ij} \geq \theta_i + \theta_j - 1, \forall i \in G', j \in G' \tag{C8}$$

$$z_{ij} \in \{0,1\}, \forall i \in G', j \in G' \tag{C9}$$

$$\theta_j \in \{0,1\}, \forall j \in G' \tag{C10}$$

$$y_{rj} \in \{0,1\}, \forall j \in G', \forall r \in R \tag{C11}$$

We generate 18 new test instances with the number of customers ranging from 10 to 1600 for both peak and non-peak scenarios. The maximum walking distance and the meeting-point separation distance are set as 1 km and 1.2 km, respectively. The problem size and computational results are shown in Table C.1 and Figure C. We can observe significantly computational time differences for larger test instances, showing the effectiveness of the layered graph structure for the hybrid metaheuristic.

Table C. Comparison of the performance of the customer-to-meeting-point assignment models with and without using the layered graph structure.

| | Without using the layered graph structure | | | | Using the layered graph structure | | | |
|---|---|---|---|---|---|---|---|---|
| $|R|$ | No. of var.* | No. of constr.* | Obj. value** | CPU(s) | No. of var.* | No of constr.* | Obj. value** | CPU(s) |
| *Peak* | | | | | | | | |
| 10 | 2.5 | 6.6 | 68.1 | 0.5 | 0.2 | 0.3 | 68.1 | 1.3 |
| 20 | 15.9 | 42.9 | 217.4 | 0.4 | 1.2 | 1.6 | 217.4 | 0.0 |
| 40 | 52.7 | 141.2 | 501.6 | 1.8 | 4.1 | 4.7 | 501.6 | 0.1 |
| 80 | 121.2 | 313.8 | 1385.5 | 5.7 | 14.5 | 16.0 | 1385.5 | 0.4 |
| 100 | 156.3 | 398.9 | 1489.4 | 7.5 | 18.4 | 16.8 | 1489.4 | 0.4 |
| 200 | 236.4 | 550.6 | 3243.6 | 12.9 | 49.4 | 37.3 | 3243.6 | 0.8 |
| 400 | 427.0 | 894.4 | 5941.8 | 21.6 | 116.1 | 65.1 | 5941.8 | 1.2 |
| 800 | 645.1 | 1140.4 | 10987.2 | 27.8 | 284.2 | 114.3 | 10987.2 | 1.9 |
| 1600 | 1079.2 | 1603.6 | 18117.4 | 31.9 | 586.7 | 176.7 | 18117.4 | 3.1 |
| *Non-peak* | | | | | | | | |
| 10 | 12.9 | 36.3 | 69.0 | 0.8 | 0.3 | 0.3 | 69.0 | 0.3 |
| 20 | 185.2 | 538.5 | 147.0 | 13.6 | 1.1 | 0.7 | 147.0 | 0.0 |
| 40 | 350.6 | 1005.6 | 314.7 | 29.2 | 3.8 | 2.0 | 314.7 | 0.0 |
| 80 | 671.6 | 1889.2 | 891.5 | 59.2 | 15.1 | 7.7 | 891.5 | 0.2 |
| 100 | 708.9 | 1967.5 | 1111.8 | 57.9 | 20.9 | 10.3 | 1111.8 | 0.2 |
| 200 | 959.1 | 2523.1 | 2815.9 | 78.3 | 73.1 | 30.1 | 2815.9 | 0.9 |
| 400 | 1135.9 | 2700.1 | 6477.9 | 84.2 | 231.1 | 80.9 | 6477.9 | 3.2 |
| 800 | 1625.8 | 3379.8 | 13895.1 | 114.8 | 617.1 | 165.9 | 13895.1 | 5.9 |
| 1600 | 2329.7 | 4036.0 | 25341.6 | 135.4 | 1407.0 | 269.8 | 25341.6 | 9.9 |

*Number of variables and number of constraints are in thousands.

** Objective function values are respect to equation (36) and (C1)

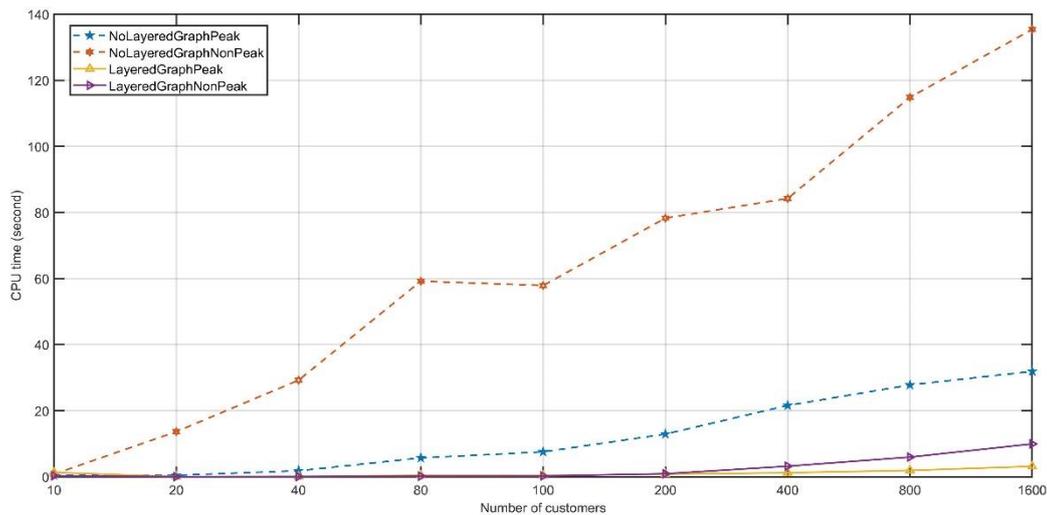

Figure C. Comparison of CPU time of the customer-to-meeting-point assignment models with and without using the layered graph structure.

## Appendix D. Computational results for solving the second-stage routing problems with charging synchronization and partial recharge

This appendix aims to assess the stand-alone effectiveness of the proposed metaheuristic for the second-stage routing problems with charging synchronization and partial recharge. Given the customer-to-meeting-point assignment obtained from the first-stage model (Eqs. (36)-(46)), the second-stage routing problem can be formulated as a MILP as follows. Let $\bar{q}_i$ denote the change in passenger load (solutions of the first stage model) at vertex $i$, $\forall i \in G' \cup D'$. By preprocessing, the set $G' \cup D'$ here contains only the assigned bus stops and transit nodes for the bus routing optimization. We denote $d_i$ the drop-off transit node corresponding to pickup bus stop node $i$. By removing the customer-to-bus stop-assignment sub-problem from Eqs (1)-(33), the second stage problem can be formulated as Eqs. (D1)-(D29). The differences from the full model are described as follows.

The objective function D1 removes total walking time of customers from Eq. (1) as the customer assignment sub-problem has been solved by the first stage model. Eq. (D6) states that if customers are picked up by a bus $k$, they must be dropped off by the same bus. Eq. (D7) updates bus passenger occupancy at bus stop and transit nodes. Eq. (D14) ensures that the bus arrives at the transit stop strictly later than it's departure from all previously visited bus stops. Eq. (D15) ensures that a customer's riding time cannot exceed their maximum ride time constraint. The other constraints are the same or equivalent to those of the full model. A preprocessing procedure is applied for removing all unused nodes and infeasible arcs to tighten the solution space of the problem. The notation is the same as the original model except the above mentioned additional variables.

We solve the above MILP problem for the 20 test instances reported in Table 4 using Gurobi solver (v11.0.0) and the metaheuristic. The results confirm the effectiveness of the metaheuristic which obtains good solutions with little computational time (31 seconds on average over the 20 test instances). Compared with the solutions obtained from Gurobi using a 4-hour computational time, the metaheuristic outperforms the commercial solver with the average gap of -11.13%. The gaps for all the test instances are almost zeros or negative.

$$\text{Min } Z = \lambda_1 \sum_{k \in K} \left( \sum_{(i,j) \in A_B} t_{ij} x_{ij}^k + \sum_{s \in S'} \tau_s^k \right) + \lambda_3 \sum_{k \in K} \sum_{i \in D'} W_i^k + \omega \sum_{i \in G'} \bar{q}_i \left( 1 - \sum_{k \in K} \sum_{j \in G' \cup D'} x_{ij}^k \right) \quad \text{(D1)}$$

$$\sum_{j \in G' \cup S' \cup \{N+1\}} x_{0j}^k = 1, \forall k \in K \quad \text{(D2)}$$

$$\sum_{i \in \{0\} \cup S' \cup D'} x_{i,N+1}^k = 1, \forall k \in K \quad \text{(D3)}$$

$$\sum_{k \in K} \sum_{i \in V_0} x_{ij}^k \leq 1, \quad \forall j \in G' \quad \text{(D4)}$$

$$\sum_{i \in V_0} x_{ij}^k - \sum_{i \in V_{N+1}} x_{ji}^k = 0, \quad \forall k \in K, j \in V \quad \text{(D5)}$$

$$\sum_{j \in V_0} x_{ji}^k \leq \sum_{j \in V} x_{jd_i}^k, \quad \forall k \in K, i \in G' \quad \text{(D6)}$$

$$x_{ij}^k = 1 \Rightarrow q_j^k = q_i^k + \bar{q}_j, \quad \forall k \in K, i \in V_0, j \in G' \cup D' \quad \text{(D7)}$$

$$0 \leq q_i^k \leq Q^k, \quad \forall k \in K, i \in V_{0,N+1} \quad \text{(D8)}$$

$$x_{ij}^k = 1 \Rightarrow B_j^k \geq B_i^k + u_i + t_{ij}, \quad \forall k \in K, i \in V_0, j \in V_{N+1} \quad \text{(D9)}$$

$$x_{sj}^k = 1 \Rightarrow B_j^k \geq B_s^k + \tau_s^k + t_{sj}, \quad \forall k \in K, s \in S', j \in \{G' \cup N+1\} \quad \text{(D10)}$$

$$x_{ij}^k = 1 \Rightarrow A_j^k = B_i^k + t_{ij} + u_i, \quad \forall k \in K, i \in G' \cup D', j \in D' \quad \text{(D11)}$$

$$p_i^k = 1 \Rightarrow W_i^k \geq B_i^k - A_i^k, \quad \forall k \in K, i \in D' \quad \text{(D12)}$$

$$p_i^k = \sum_{j \in V} x_{ji}^k, \quad \forall k \in K, i \in D' \quad \text{(D13)}$$

$$\sum_{j\in G'\cup D'} x_{ji}^k = 1 \Rightarrow A_{d_i}^k \geq B_i^k + t_{id_i} + u_i, \quad \forall k \in K, i \in G' \tag{D14}$$

$$\sum_{j\in G'\cup D'} x_{ji}^k = 1 \Rightarrow A_{d_i}^k - (B_i^k + u_i) \leq L_i, \quad \forall k \in K, i \in G' \tag{D15}$$

$$e_i \leq B_i^k \leq l_i, \quad \forall k \in K, i \in V \tag{D16}$$

$$E_0^k = E_{init}^k, \quad \forall k \in K \tag{D17}$$

$$E_{min}^k \leq E_i^k \leq E_{max}^k, \quad \forall k \in K, i \in V \tag{D18}$$

$$x_{ij}^k = 1 \Rightarrow E_j^k = E_i^k - \beta^k c_{ij}, \quad \forall k \in K, i \in V_0\setminus S', j \in V_{N+1} \tag{D19}$$

$$x_{ij}^k = 1 \Rightarrow E_j^k = E_s^k + \alpha_s \tau_s^k - \beta^k c_{sj}, \quad \forall k \in K, s \in S', j \in \{G' \cup N+1\} \tag{D20}$$

$$v_s = \sum_{k\in K}\sum_{j\in G'\cup N+1} x_{sj}^k, s \in S' \tag{D21}$$

$$v_h \leq v_l, \forall h,l \in S'_o, o \in S, h < l \tag{D22}$$

$$v_h + v_l = 2 \Rightarrow \sum_{k\in K} B_h^k \geq \sum_{k\in K} B_l^k + \sum_{k\in K} \tau_l^k, \forall h,l \in S'_o, o \in S, h < l \tag{D23}$$

$$\tau_s^k + B_s^k \leq M_2 \sum_{j\in G'\cup N+1} x_{sj}^k, \forall s \in S', k \in K \tag{D24}$$

$$v_s \leq 1, s \in S' \tag{D25}$$

$$x_{ij}^k \in \{0,1\}, \forall k \in K, i,j \in V_{0,N+1} \tag{D26}$$

$$\tau_s^k \geq 0, v_s \in \{0,1\}, \forall k \in K, s \in S' \tag{D27}$$

$$A_i^k \geq 0, B_i^k \geq 0, \forall k \in K, i \in V_{0,N+1} \tag{D28}$$

$$p_i^k \in \{0,1\}, W_i^k \geq 0, \forall k \in K, i \in D' \tag{D29}$$

**Table D. Computational results for solving the second-stage routing problems for the test instances.**

| Instance | MILP | | | | Hybrid metaheuristic | | | | |
|---|---|---|---|---|---|---|---|---|---|
| | Best known solution (BKS) | Gap to the lower bound | CT* (min) | CPU (s) | Best** obj. | Gap*** (avg.) | Gap*** (best) | CT* (min) | CPU (s) |
| c10op | 64.38 | 2.95% | 4.9 | 14400 | 64.38 | 0.00% | 0.00% | 4.9 | 4 |
| c20op | 131.70 | 30.48% | 18.2 | 14400 | 131.70 | 0.06% | 0.00% | 18.2 | 5 |
| c30op | 161.60 | 40.76% | 16.7 | 14400 | 161.49 | 0.25% | -0.06% | 16.5 | 15 |
| c40op | 231.01 | 49.35% | 30.2 | 14400 | 231.15 | 0.42% | 0.06% | 30.4 | 13 |
| c50op | 290.10 | 50.64% | 43.3 | 14400 | 289.49 | 0.03% | -0.21% | 42.7 | 14 |
| c60op | 457.36 | 67.87% | 28.5 | 14400 | 307.37 | -32.58% | -32.79% | 26.3 | 53 |
| c70op | 467.45 | 67.70% | 37.6 | 14400 | 337.21 | -27.51% | -27.86% | 32.5 | 70 |
| c80op | 508.27 | 68.22% | 63.1 | 14400 | 385.53 | -24.01% | -24.15% | 42.5 | 46 |
| c90op | 441.20 | 61.45% | 41.1 | 14400 | 403.57 | -7.96% | -8.53% | 31.7 | 108 |
| c100op | 426.12 | 63.47% | 37.2 | 14400 | 385.84 | -8.69% | -9.45% | 27.4 | 78 |
| Average | 317.92 | 50.29% | 32.07 | 14400 | 269.77 | -10.00% | -10.30% | 27.3 | 41 |
| c10p | 64.57 | 10.62% | 4.9 | 14400 | 64.57 | 0.00% | 0.00% | 4.9 | 5 |
| c20p | 215.95 | 65.23% | 11.3 | 14400 | 215.95 | 0.04% | 0.00% | 11.3 | 6 |
| c30p | 158.47 | 58.56% | 8.2 | 14400 | 158.46 | 0.01% | -0.01% | 8.2 | 12 |
| c40p | 181.57 | 58.10% | 13.1 | 14400 | 181.33 | -0.09% | -0.13% | 12.9 | 13 |
| c50p | 288.90 | 71.25% | 26.5 | 14400 | 288.90 | 0.08% | 0.00% | 26.5 | 12 |
| c60p | 193.80 | 52.62% | 4.1 | 14400 | 193.23 | 0.31% | -0.30% | 4.0 | 29 |
| c70p | 232.16 | 58.90% | 11.4 | 14400 | 231.90 | 0.03% | -0.11% | 11.6 | 33 |
| c80p | 217.95 | 57.70% | 9.8 | 14400 | 217.31 | -0.08% | -0.29% | 8.9 | 26 |
| c90p | 787.15 | 86.97% | 15.4 | 14400 | 339.96 | -56.70% | -56.81% | 19.3 | 34 |
| c100p | 946.01 | 89.49% | 12.5 | 14400 | 317.90 | -66.27% | -66.40% | 14.2 | 37 |
| Average | 328.65 | 60.94% | 11.72 | 14400 | 220.95 | -12.27% | -12.40% | 12.2 | 21 |
| Overall average | 323.29 | 55.62% | 21.90 | 14400 | 245.36 | -11.13% | -11.35% | 19.7 | 31 |

\* CT: charging time
\*\* Based on the average of 5 runs. The gaps are compared to the best solutions found by Gurobi (the first column).
   Charging time of the obtained solutions is measured in minutes. CPU time is measured in seconds.